\documentclass[12 pt]{amsart}

\usepackage[utf8]{inputenc}
\usepackage[all]{xy}
\usepackage{bm}
\usepackage{enumitem}
\usepackage{colonequals,amsmath,amssymb,amsthm,mathtools}
\usepackage{mathrsfs}
\usepackage[ruled,vlined]{algorithm2e}
\usepackage{float}
\usepackage{tikz,tikz-cd,graphicx}
\graphicspath{ {./images/} }
\usepackage{caption}
\usepackage{todonotes}
\usepackage[colorlinks=true]{hyperref}
\hypersetup{colorlinks=true,  linkcolor=blue,  filecolor=blue,  citecolor=red, urlcolor=blue}
\usepackage[left=30mm, right=30mm, top=25mm, bottom=25mm]{geometry}
\usepackage[capitalize,noabbrev]{cleveref}
\usepackage{pdfpages}
\usepackage{microtype}
\usepackage[justification=centering]{caption}
\usepackage{indentfirst}

\newtheorem{theorem}{Theorem}

\newtheorem{corollary}[theorem]{Corollary}

\theoremstyle{definition}
\newtheorem{definition}[theorem]{Definition}

\newtheorem{remark}[theorem]{Remark}

\newtheorem{example}[theorem]{Example}
\newtheorem{assumption}[theorem]{Assumption}

\definecolor{ggreen}{RGB}{0, 94, 25}

\newcommand{\an}{{\operatorname{an}}}

\newcommand{\GL}{{\operatorname{GL}}}

\newcommand{\trop}{{\operatorname{trop}}}

\newcommand{\Spec}{{\operatorname{Spec}}}

\newcommand{\SL}{{\operatorname{SL}}}

\newcommand*\angles[1]{\langle #1 \rangle}

\newcommand*\Z{\mathbb{Z}}

\newcommand*\R{\mathbb{R}}

\newcommand*\A{\mathbb{A}}
\renewcommand*\P{\mathbb{P}}
\newcommand*\T{\mathbb{T}}

\newcommand*\mfrak{\mathfrak{m}}

\newcommand*\wb{\overline} 
\newcommand*\varhrulefill[1][3pt]{\leavevmode\leaders\hrule height#1\hfill\kern0pt}

\newcommand*{\DashedArrow}[1][]{\mathbin{\tikz [baseline=-0.25ex,-latex, dashed,#1] \draw [#1] (0pt,0.5ex) -- (1.3em,0.5ex);}}

\definecolor{lightyellow}{RGB}{255, 255, 197}

\title[Tropical invariants for binary quintics and reduction  of Picard curves]{Tropical invariants for binary quintics and reduction types of Picard curves}

\author{Paul Alexander Helminck}
\address{Department of Mathematical Sciences Durham University,	Upper Mountjoy Campus,	Stockton Road,	Durham DH1 3LE.}
\email{paul.a.helminck@durham.ac.uk}

\author{Yassine El Maazouz}
\address{U.C. Berkeley, Department of Statistics, 335 Evans Hall \# 3860, Berkeley CA 94720, U.S.A.}
\email{yassine.el-maazouz@berkeley.edu}

\author{Enis Kaya}
\address{Department of Mathematics, KU Leuven, Celestijnenlaan 200B, 3001 Heverlee, Belgium}
\email{enis.kaya@kuleuven.be}

\begin{document}

	\maketitle
	
\begin{abstract}
    In this paper, we express the reduction types of Picard curves in terms of tropical invariants associated to binary quintics. We also give a general framework for tropical invariants associated to group actions on arbitrary varieties. The problem of finding tropical invariants for binary forms fits in this general framework by mapping the space of binary forms to symmetrized versions of the Deligne--Mumford compactification~$\overline{M}_{0,n}$.
\end{abstract}

\section{Introduction}\label{Sec:Intro}

Invariant theory studies quantities in geometry that are invariant under group actions. This theory sparked many developments in commutative algebra, leading to the Hilbert basis theorem and many other results. In this paper, we study invariants of binary forms $f(x,z)=a_{0}x^{n} + a_{1} x^{n-1}z + \dots + a_{n}z^{n}$ defined over an algebraically closed field $K$ of characteristic $0$. The group in question is $\GL_{2}(K)$ and it acts on these binary forms through \emph{M\"{o}bius transformations}; that is
\[
    f^{\sigma}(x,z) = f(ax+bz,cx+dz), \quad \text{for } \sigma=\begin{bmatrix}
    a & b \\
    c & d
    \end{bmatrix} \in \GL_2(K).
\]
By mapping the coefficients $a_i$ of $f$ to those of $f^{\sigma}$, this then also gives an action of $\GL_{2}(K)$ on the ring $A \coloneqq K[a_{0},\dots,a_{n}]$. The \emph{invariants} for this action are homogeneous polynomials $H \in A$ such that for any $\sigma \in \GL_2(K)$ one has $H^{\sigma}=\mathrm{det}(\sigma)^{k}H$ for some $k \in \Z$. Similarly, we say that a homogeneous polynomial $H \in A$ is $\SL_{2}(K)$-invariant if $H^{\sigma} = H$ for all $\sigma \in \SL_{2}(K)$. Note that these two types of invariants are the same. We will from now on restrict our attention to the $G \coloneqq \SL_{2}(K)$ action.
By well-known results from commutative algebra, the invariant polynomials form a finitely generated subring $A^G$ of $A$, called the ring of invariants. There are algorithms that can explicitly calculate the generators of this ring and a full list of generators is known for all binary forms of degree $n \leq 10$. These generators of the ring of invariants satisfy the pleasant property that two separable binary forms lie in the same $\GL_{2}(K)$-orbit if and only if their invariants define the same point in the weighted projective space associated to the generators, where the weights are the degrees of the generators.
This equivalence gives a strong connection between geometry on the one hand and algebra on the other.

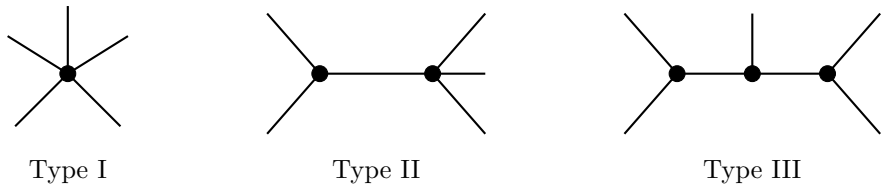
\begin{figure}[H]
    	    \centering
            	    \begin{tikzpicture}
                        \filldraw (0.9,0) circle (3pt);
                        \draw [thick] (0.9,0) -- (0.9,0.9);
                        \draw [thick] (0.9,0) -- (1.7,0.5);
                        \draw [thick] (0.9,0) -- (0.1,0.5);
                        \draw [thick] (0.9,0) -- (1.6,-0.7);
                        \draw [thick] (0.9,0) -- (0.2,-0.7);
                        
                        \filldraw (0.9,-1.3) node{Type~$\mathrm{I}$};
                        
                        \qquad
                        
                        \filldraw (4.25,0) circle (3pt);
                        \filldraw (5.75,0) circle (3pt);
                        \draw [thick] (4.25,0) -- (5.75,0);
                        \draw [thick] (4.25,0) -- (3.55,0.8);
                        \draw [thick] (4.25,0) -- (3.55,-0.8);
                        \draw [thick] (5.75,0) -- (6.45,0.8);
                        \draw [thick] (5.75,0) -- (6.45,0);
                        \draw [thick] (5.75,0) -- (6.45,-0.8);
                
                        \filldraw (5,-1.3) node{Type~$\mathrm{II}$};
                
                         \qquad
                         
                        \filldraw (9,0) circle (3pt);
                        \filldraw (10,0) circle (3pt);
                        \filldraw (11,0) circle (3pt);
                        \draw [thick] (9,0) -- (11,0);
                        \draw [thick] (9,0) -- (8.3,0.8);
                        \draw [thick] (9,0) -- (8.3,-0.8);
                        \draw [thick] (10,0) -- (10,0.8);
                        \draw [thick] (11,0) -- (11.7,0.8);
                        \draw [thick] (11,0) -- (11.7,-0.8);
                
                        \filldraw (10,-1.3) node{Type~$\mathrm{III}$};
                
                        \end{tikzpicture}
    	    \caption{The three types of unmarked phylogenetic trees with $5$ leaves}
    	    \label{fig:tropicalTypes}
    	\end{figure}

We now turn to the non-archimedean side of this story and consider a complete non-archimedean algebraically closed field $K$ of characteristic zero with non-trivial valuation $v\colon K^{*} \to \mathbb{R}$. Let $f(x,z)$ be a separable binary form of degree $n$ over $K$. The zeroes of this form give a canonical metric tree on $n$ leaves by connecting the corresponding points in the Berkovich analytification of $\mathbb{P}^{1}$. There are only finitely many phylogenetic types for any given degree $n$ (see \cref{fig:tropicalTypes} for the case $n=5$), and the possible types give a partition of the space of all non-archimedean binary forms of degree $n$. Since the tree is invariant under projective isomorphisms, this also partitions the space of all invariants. A natural question now arises: what are the equations for these partitions?

There are two instances where we know the equations: $n=4$ and $n=6$. For $n=4$ there are two tree types: a trivial one with four leaves connected to a single vertex and a non-trivial one. We can distinguish between these two using the valuation of the $j$-invariant of the quartic. Namely, a quartic has trivial tree type if and only if $v(j)\geq{0}$, see \cite[Chapter VII]{SilvermanI}. For $n=6$, there are seven tree types and one can distinguish between them using the valuations of Igusa invariants, see \cite{Liu_CourbesStablesGenre2} and \cite{helminck2021Igusa}. Our main goal in this paper is to fill up this gap and find the invariants for \emph{binary quintics}; that is, for $n=5$.

For quintics, there are three tree types; see Figure \ref{fig:tropicalTypes}. We wish to distinguish between these tree types using the valuations of suitable invariants. To that end, we start with a set of generators $I_{4},I_{8},I_{12},I_{18}$ for the ring of invariants together with the discriminant $\Delta$. The valuations of these invariants are not sufficient to determine the tree type of the quintic, so we introduce a new invariant $H$, giving the set $S=\{I_{4},I_{8},I_{12},I_{18},\Delta,H\}$. In our first theorem, we show that the valuations of these invariants determine the tree type of a quintic. We call the valuations of these invariants the \emph{tropical invariants} of the quintic. For technical reasons, we assume for the remainder of this section that the residue characteristic $p$ of $K$ is not equal to $2,3,11$. Note however that Remark~\ref{rem:Thm1p=11} explains how to deal with the case where $p=11$.

    \begin{theorem}[{\bf Tree types of binary quintics}] \label{thm:tropicalTypes}
        Let $f$ be a separable binary quintic over $K$. Then, the tree type of $f$ is determined by the tropical invariants as follows: 
    	\begin{enumerate}[label=(\Roman*)]
        	\item \label{MT1:type1} 
          	        The tree is of Type~$\mathrm{I}$ if and only if $8v(I)-\mathrm{deg}(I)v(\Delta)\geq{0}$ for all $I\in{S}$.
    	    \item \label{MT1:type2}
    	            The tree is of Type~$\mathrm{II}$ if and only if $v(\Delta)-2v(I_{4})>0$ or $9v(\Delta)-4v(I_{18})>0$, and $12v(I)-\mathrm{deg}(I)v(H)\geq{0}$ for all $I\in{S}$.
    	    \item \label{MT1:type3}
    	            The tree is of Type~$\mathrm{III}$ if and only if $v(\Delta)-2v(I_{4})>0$ and $v(H)-3v(I_{4})>0$.  
	    \end{enumerate}
    \end{theorem}
    
\begin{figure}[ht]
    	    \centering
            	    \begin{tikzpicture}
                        \filldraw (0.9,0) circle (3pt);
                        \draw [thick] (0.9,0) -- (0.9,0.9);
                        \draw [thick] (0.9,0) -- (1.7,0.5);
                        \draw [thick] (0.9,0) -- (0.1,0.5);
                        \draw [thick] (0.9,0) -- (1.6,-0.7);
                        \draw [thick] (0.9,0) -- (0.2,-0.7);
                        \filldraw[blue] (0.9,0.9) circle (2pt);
                        \filldraw (0.9,1.3) node{$\infty$};
                        
                        \filldraw (0.9,-1.3) node{Type~$\mathrm{I}$};
                        
                        \qquad
                        
                        \filldraw (4.25,1.5) circle (3pt);
                        \filldraw (5.75,1.5) circle (3pt);
                        \draw [thick] (4.25,1.5) -- (5.75,1.5);
                        \draw [thick] (4.25,1.5) -- (3.55,2.3);
                        \draw [thick] (4.25,1.5) -- (3.55,0.7);
                        \draw [thick] (5.75,1.5) -- (6.45,2.3);
                        \draw [thick] (5.75,1.5) -- (6.45,1.5);
                        \draw [thick] (5.75,1.5) -- (6.45,0.7);
                        \filldraw[blue] (3.55,2.3) circle (2pt);
                        \filldraw (3.15,2.55) node{$\infty$};
                
                        \filldraw (5,0.3) node{Type~$\mathrm{II}.1$};
                        
                        \filldraw (4.25,-1.5) circle (3pt);
                        \filldraw (5.75,-1.5) circle (3pt);
                        \draw [thick] (4.25,-1.5) -- (5.75,-1.5);
                        \draw [thick] (4.25,-1.5) -- (3.55,-0.7);
                        \draw [thick] (4.25,-1.5) -- (3.55,-2.3);
                        \draw [thick] (5.75,-1.5) -- (6.45,-0.7);
                        \draw [thick] (5.75,-1.5) -- (6.45,-1.5);
                        \draw [thick] (5.75,-1.5) -- (6.45,-2.3);
                        \filldraw[blue] (6.45,-0.7) circle (2pt);
                        \filldraw (6.85,-0.45) node{$\infty$};
                
                        \filldraw (5,-2.7) node{Type~$\mathrm{II}.2$};
                
                         \qquad
                         
                        \filldraw (9,1.5) circle (3pt);
                        \filldraw (10,1.5) circle (3pt);
                        \filldraw (11,1.5) circle (3pt);
                        \draw [thick] (9,1.5) -- (11,1.5);
                        \draw [thick] (9,1.5) -- (8.3,2.3);
                        \draw [thick] (9,1.5) -- (8.3,0.7);
                        \draw [thick] (10,1.5) -- (10,2.3);
                        \draw [thick] (11,1.5) -- (11.7,2.3);
                        \draw [thick] (11,1.5) -- (11.7,0.7);
                        \filldraw[blue] (8.3,2.3) circle (2pt);
                        \filldraw (7.9,2.55) node{$\infty$};
                
                        \filldraw (10,0.3) node{Type~$\mathrm{III}.1$};
                
                        \filldraw (9,-1.5) circle (3pt);
                        \filldraw (10,-1.5) circle (3pt);
                        \filldraw (11,-1.5) circle (3pt);
                        \draw [thick] (9,-1.5) -- (11,-1.5);
                        \draw [thick] (9,-1.5) -- (8.3,-0.7);
                        \draw [thick] (9,-1.5) -- (8.3,-2.3);
                        \draw [thick] (10,-1.5) -- (10,-0.7);
                        \draw [thick] (11,-1.5) -- (11.7,-0.7);
                        \draw [thick] (11,-1.5) -- (11.7,-2.3);
                        \filldraw[blue] (10,-0.7) circle (2pt);
                        \filldraw (10,-0.3) node{$\infty$};
                
                        \filldraw (10,-2.7) node{Type~$\mathrm{III}.2$};
                
                        \end{tikzpicture}
    	    \caption{Tree types of binary $(4,1)$-forms}
    	    \label{fig:fourOneTreeTypes}
    	\end{figure}
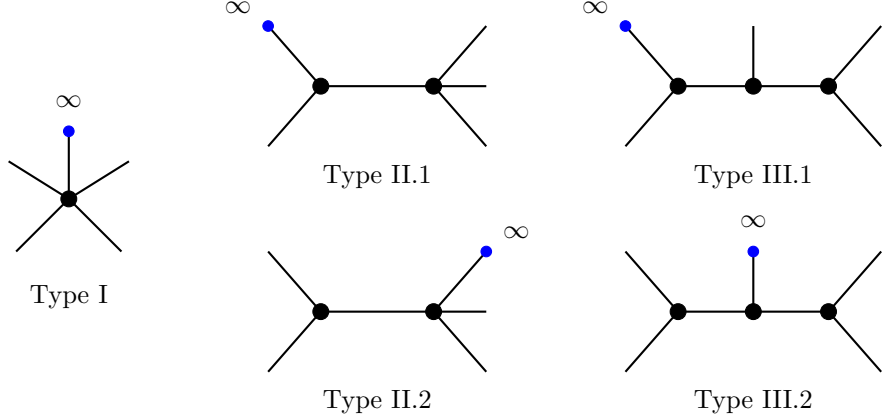

Finding the equations for $n=4$ and $n=6$ is partially motivated by applications to reduction types of hyperelliptic curves. For $n=5$, the motivation comes from Picard curves $X$, which are smooth plane quartics of the form 
\begin{equation*}
    y^3\ell(x,z)=q(x,z),
\end{equation*}
where $\deg(\ell(x,z))=1$ and $\deg(q(x,z))=4$. We are interested in obtaining the \emph{minimal skeleton} of the Berkovich analytification of such a curve, which codifies the different possible semistable models for $X$, see \cite[Section~4.16]{Baker_Payne_Rabinoff_13}. The results in \cite{helminckInvariantsTreesSuper} show that this skeleton can be recovered from the \emph{marked tree type} of the  quintic $f(x,z) = \ell(x,z)\cdot {q(x,z)}$. More precisely, this quintic has five distinct roots, giving a metric tree with five leaves, and the root of $\ell(x,z)$ gives the marking. We can assume by a projective transformation that this marked point is $\infty$. There are exactly five marked tree types, see Figure \ref{fig:fourOneTreeTypes}, giving rise to five reduction types of Picard curves. In terms of invariant theory, the natural object to consider here is the binary $(4,1)$-form $(q(x,z),\ell(x,z))$. These binary $(4,1)$-forms similarly have a finitely generated ring of invariants and in our second theorem we show that we can find a set of invariants for $(4,1)$-forms that distinguish between the five marked tree types. This set consists of the set $S$ from Theorem~\ref{thm:tropicalTypes} together with a new set $S'=\{j_{2},j_{3},j_{5},j_{6},j_{9}\}$ of $(4,1)$-invariants. As above, we call the valuations of these invariants the \emph{tropical invariants} of the $(4,1)$-form. 

     \begin{theorem}[{\bf Tree types of $(4,1)$-forms}]  \label{thm:reductionTypes} 
            Let $(q,\ell)$ be a $(4,1)$-form over $K$ such that the associated binary quintic $f = q\cdot\ell$ is separable. The tree type of $(q,\ell)$ is determined by the tropical invariants as follows: 
            \begin{enumerate}[label=(\Roman*)]
                \item If $f$ has tree Type~$\mathrm{I}$, then $(q,\ell)$ also has Type~$\mathrm{I}$.
            
                \item If $f$ has tree Type~$\mathrm{II}$, then $(q,\ell)$ has Type~$\mathrm{II}.1$ (resp. Type~$\mathrm{II}.2$) if and only if the quantity $5 v(j_{2})- 2v(j_{5})$ is strictly positive (resp. zero).
                
                \item If $f$ has tree Type~$\mathrm{III}$, then $(q,\ell)$ has Type~$\mathrm{III}.1$ (resp. Type~$\mathrm{III}.2$) if and only if the quantity $5 v(j_{2}) -  2v(j_{5})$ is strictly positive (resp. zero).
            \end{enumerate}
    \end{theorem}

To determine the skeleton of a Picard curve, we also need to know the associated weights and edge lengths. The weights of the skeleton are completely determined by \cref{thm:reductionTypes}, but the lengths are not. In our third theorem, we give formulas for the edge lengths of a $(4,1)$-form in terms of its tropical invariants. For trees of Type~$\mathrm{II}$ and Type~$\mathrm{III}.2$, we are able to give these in terms of invariants of quintics. For trees of Type~$\mathrm{III}.1$, we express the marked edge lengths in terms of $(4,1)$-invariants. This difference is quite natural, since there is a natural symmetry on trees of Type~$\mathrm{III}.2$.

\begin{theorem}[{\bf Edge lengths}] \label{thm:edgeLengths}
    The edge lengths of the trees in \cref{thm:reductionTypes} are given by the tropical invariants as follows:
        \begin{enumerate}[label=(\Roman*)]
            \item If $f$ has Type~$\mathrm{I}$, then there are no non-trivial edges.
        
            \item If $f$ has Type~$\mathrm{II}$, then there is only one edge $e_1$. Its length, both in the cases of unmarked and marked trees, is given by
            \[
            L(e_1) = \max\left(\frac{1}{2}(v(\Delta) - 2v(I_4)),\  \frac{1}{3}(2 v(\Delta) - v(I_{18})) \right).
            \]
            
            \item If $f$ has Type~$\mathrm{III}$, then there are two edges $e_{1}$ and $e_{2}$. Assume the length of $e_1$ is less than or equal to that of $e_2$. The unmarked edge lengths are given by
            \begin{equation}\label{eq:EdgeLengthsTypeIII}
            \begin{aligned}
                L(e_{1})&=\mathrm{min}\left( \frac{1}{2}\left(v(I_{18}) - \frac{9}{2} v(I_{4})\right),\frac{1}{4}\left(v(\Delta)-2v(I_{4})\right) \right),\\
                L(e_{2})&=v(\Delta)-2v(I_{4})-2L(e_{1}).
            \end{aligned}
            \end{equation}
            If $(q,\ell)$ has Type~$\mathrm{III}.1$, then we write $e_{1}$ for the edge adjacent to the marked point and $e_{2}$ for the other edge. The edge lengths are then given by
            \begin{align*}
                L(e_{1}) &= \frac{1}{10}(5v(j_{2})-2v(j_{5})),\\
                L(e_{2}) &= \frac{1}{2}(v(\Delta)-2v(I_{4}))-L(e_{1}).
            \end{align*}
            For trees of Type~$\mathrm{III}.2$, they are as
            in \eqref{eq:EdgeLengthsTypeIII}.
        \end{enumerate}
\end{theorem} 

    Combining these theorems, we then immediately obtain a description of the reduction types of Picard curves in terms of quintic and $(4,1)$-invariants.
    
        \begin{figure}[ht]
    	    \centering
            	    \begin{tikzpicture}
                        \filldraw (0.9,0) circle (3pt);
                        \draw [thick] (0.9,0) -- (0.9,0.9);
                        \draw [thick] (0.9,0) -- (1.7,0.5);
                        \draw [thick] (0.9,0) -- (0.1,0.5);
                        \draw [thick] (0.9,0) -- (1.6,-0.7);
                        \draw [thick] (0.9,0) -- (0.2,-0.7);
                        \filldraw[blue] (0.9,0.9) circle (2pt);
                        \filldraw (0.9,1.3) node{$\infty$};
                        \filldraw (0.9,-0.5) node{$3$};
                        
                        \filldraw (0.9,-1.3) node{Type~$\mathrm{I}$};
                        
                        \qquad
                        
                        \filldraw (4.25,1.5) circle (3pt);
                        \filldraw (5.75,1.5) circle (3pt);
                        \draw [thick] (4.25,1.5) -- (5.75,1.5);
                        \draw [thick] (4.25,1.5) -- (3.55,2.3);
                        \draw [thick] (4.25,1.5) -- (3.55,0.7);
                        \draw [thick] (5.75,1.5) -- (6.45,2.3);
                        \draw [thick] (5.75,1.5) -- (6.45,1.5);
                        \draw [thick] (5.75,1.5) -- (6.45,0.7);
                        \filldraw[blue] (3.55,2.3) circle (2pt);
                        \filldraw (3.15,2.55) node{$\infty$};
                        \draw[thick] (5,1.5) ellipse (0.75 and 0.42);
                        \filldraw (5.75,1) node{$1$};
                
                        \filldraw (5,0.3) node{Type~$\mathrm{II}.1$};
                        
                        \filldraw (4.25,-1.5) circle (3pt);
                        \filldraw (5.75,-1.5) circle (3pt);
                        \draw [thick] (4.25,-1.5) -- (5.75,-1.5);
                        \draw [thick] (4.25,-1.5) -- (3.55,-0.7);
                        \draw [thick] (4.25,-1.5) -- (3.55,-2.3);
                        \draw [thick] (5.75,-1.5) -- (6.45,-0.7);
                        \draw [thick] (5.75,-1.5) -- (6.45,-1.5);
                        \draw [thick] (5.75,-1.5) -- (6.45,-2.3);
                        \filldraw[blue] (6.45,-0.7) circle (2pt);
                        \filldraw (6.85,-0.45) node{$\infty$};
                        \filldraw (4.25,-2) node{$1$};
                        \filldraw (5.75,-2) node{$2$};
                
                        \filldraw (5,-2.7) node{Type~$\mathrm{II}.2$};
                
                         \qquad
                         
                        \filldraw (9,1.5) circle (3pt);
                        \filldraw (10,1.5) circle (3pt);
                        \filldraw (11,1.5) circle (3pt);
                        \draw [thick] (9,1.5) -- (11,1.5);
                        \draw [thick] (9,1.5) -- (8.3,2.3);
                        \draw [thick] (9,1.5) -- (8.3,0.7);
                        \draw [thick] (10,1.5) -- (10,2.3);
                        \draw [thick] (11,1.5) -- (11.7,2.3);
                        \draw [thick] (11,1.5) -- (11.7,0.7);
                        \filldraw[blue] (8.3,2.3) circle (2pt);
                        \filldraw (7.9,2.55) node{$\infty$};
                        \draw[thick] (9.5,1.5) ellipse (0.50 and 0.32);
                        \filldraw (11,1) node{$1$};
                
                        \filldraw (10,0.3) node{Type~$\mathrm{III}.1$};
                
                        \filldraw (9,-1.5) circle (3pt);
                        \filldraw (10,-1.5) circle (3pt);
                        \filldraw (11,-1.5) circle (3pt);
                        \draw [thick] (9,-1.5) -- (11,-1.5);
                        \draw [thick] (9,-1.5) -- (8.3,-0.7);
                        \draw [thick] (9,-1.5) -- (8.3,-2.3);
                        \draw [thick] (10,-1.5) -- (10,-0.7);
                        \draw [thick] (11,-1.5) -- (11.7,-0.7);
                        \draw [thick] (11,-1.5) -- (11.7,-2.3);
                        \filldraw[blue] (10,-0.7) circle (2pt);
                        \filldraw (10,-0.3) node{$\infty$};
                        \filldraw (9,-2) node{$1$};
                        \filldraw (10,-2) node{$1$};
                        \filldraw (11,-2) node{$1$};
                
                        \filldraw (10,-2.7) node{Type~$\mathrm{III}.2$};
                
                        \end{tikzpicture}
    	    \caption{Reduction types of Picard curves}
    	    \label{fig:ReductionTypesOfPicardCurves}
    	\end{figure}
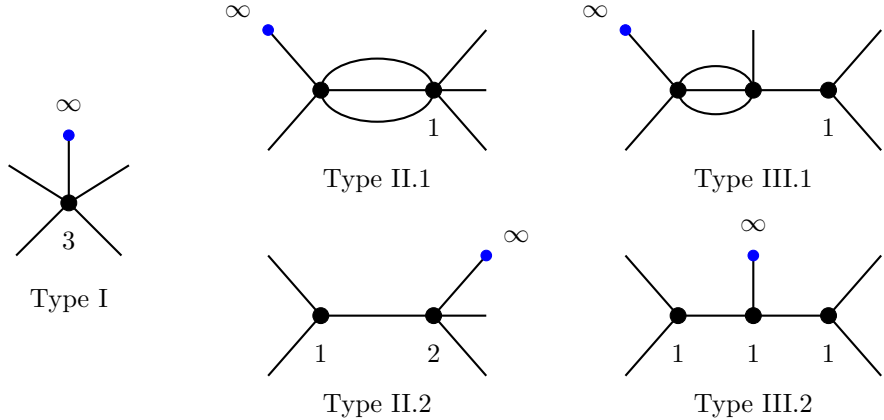

        \begin{corollary}[{\bf Reduction types of Picard curves}]\label{cor:reductionTypesOfPicardCurves}
            The reduction type of the Picard curve $y^3 \ell(x,z) = q(x,z)$ is completely determined by the tropical invariants of the $(4,1)$-form $(q,\ell)$. The tree types given in \cref{fig:fourOneTreeTypes} correspond to the reduction types in \cref{fig:ReductionTypesOfPicardCurves}.
        \end{corollary}

        We thus obtain a description of the moduli space of tropical Picard curves in terms of invariants of quintics and $(4,1)$-forms. The work in \cite{Clery_vanderGeer_21} shows that these invariants in fact give rise to Picard modular forms. This is completely analogous to the cases $n=4$ and $n=6$ mentioned above. For instance, for $n=4$ we have the invariants $c_{4}$ and $c_{6}$ corresponding to Eisenstein modular forms for $\SL_{2}(\mathbb{Z})$. The $j$-invariant is a rational function in terms of these modular forms and the classical criterion an elliptic curve $E$ has (potential) good reduction if and only if its $j$-invariant is non-negative expresses the tropical moduli space in terms of these modular forms. For $n=6$ the Igusa invariants similarly give rise to modular forms (see \cite[Section~6]{vanderGeer2021Siegel}) and the criteria given in \cite{Liu_CourbesStablesGenre2} and \cite{helminck2021Igusa} again express the tropical moduli space in terms of modular forms. We view our results as extensions of those for Picard curves.

        \subsection{Comparison with existing literature}
        
        The results obtained in this paper can be seen as a natural continuation of \cite{Liu_CourbesStablesGenre2} and \cite{helminck2021Igusa}. In the first, criteria for the seven reduction types of curves of genus two were given in terms of the Igusa invariants. In \cite{helminck2021Igusa}, this result was extended to arbitrary complete non-archimedean fields and an easier proof was given. This paper in turn was based on \cite{helminckInvariantsTreesSuper}, where skeleta of general superelliptic curves are studied. In the latter, it was shown that one can recover the skeleton from tropicalizations of certain functions in the coefficients of $f(x)$. The key difference between this paper and the latter is that the functions we give here are projective invariants of binary forms. This can be used to interpret the criteria in terms of Picard modular forms, giving a stronger connection to various moduli spaces in the literature. 
        
        We also define a general notion of a set of tropical invariants using tools from non-archimedean geometry. A notion that seems distantly related to this one appears in \cite{KavehManon2019}, where spherical varieties and invariant valuations are studied. We do not restrict ourselves to spherical varieties, as we can phrase everything in terms of $G$-invariant subsets of the analytification of an algebraic variety for some group $G$.
        
        From a geometric point of view, our paper fits into the literature as follows. Suppose that we have a group action on a variety that admits a geometric quotient. This quotient admits many possible compactifications and for every compactification we obtain a natural definition of a tropical invariant. For instance, for separable binary forms we can compactify the space using either stable binary forms (see \cite{MumfordFogartyKirwan1994}) or the symmetrized Deligne--Mumford compactification $\overline{M}_{0,n}/S_{n}$. We are mostly interested in the latter, since this has direct applications to reduction types of Picard curves. That is, if we write $\overline{\mathcal{N}}$ for the space of admissible $\mathbb{Z}/3\mathbb{Z}$-coverings with ramification signature $(4,1)$, then there is a natural map $\overline{\mathcal{N}}\to \overline{M}_{0,5} / (S_{4} \times S_{1})$ sending a covering to its branch locus and this map respects the boundary loci. Recent work by Cl{\'e}ry and van der Geer \cite{Clery_vanderGeer_21} shows that this map can be used to connect $(4,1)$-invariants to Picard modular forms, mirroring the classical case of elliptic curves. It would be interesting to see how tropical invariants of binary forms are connected to other moduli spaces.

        In classical invariant theory, many important concepts such as Hilbert series and generators of invariant rings can be computed explicitly using Gr\"{o}bner bases. In tropical geometry, Gr\"{o}bner bases have similarly led to many explicit algorithms as well. This analogy suggests the existence of an algorithm to compute tropical invariants using Gr\"{o}bner-theoretic methods. We leave an in-depth study of this topic for future work.

        \subsection{Outline}
        
	        This paper is organized as follows. In \cref{Sec:Background}, we review some background on invariant theory for binary forms and Picard curves. \cref{Sec:tropicalInvariants} introduces and discusses  the notion of tropical invariants. Finally, we prove our main results, and in particular discuss the edge lengths (or thickness of singular points), in \cref{Sec:ProofMainResults}. Finally we note that several results in this article were found or proved by symbolic computations. The codes and computations are made available at
            \begin{equation}\label{eq:link_for_code}
            \text{\url{https://mathrepo.mis.mpg.de/TropicalInvariantsPicardCurves/index.html}}.
            \end{equation}
            Our code is implemented in SageMath \cite{sagemath}.

	   \subsection{Notation}
	
        	Throughout this paper, unless explicitly stated otherwise, we will use the following notation:
        	
        	\vspace{1mm}
        	
        	\begin{center}
            	\begin{tabular}{rcl}
                		$K$ &--& Complete, non-archimedean and algebraically closed \\
                		    &  & valued field of characteristic $0$.\\
                		$v$ &--& Non-trivial valuation on $K$.\\
                		$R$ &--& The valuation ring of $K$.\\
                		$\mfrak$ &--& The unique maximal ideal of $R$.\\
                		$k$ &--& The residue field $R/\mfrak$. \\
                		$p$ &--& The residue characteristic of $K$; that is $p = \mathrm{char}(k)$.\\
                		$\A^n$ &--& Affine space over $K$ of dimension $n$. \\
                		$\P^n$ &--& Projective space over $K$ of dimension $n$. \\
                		$\wb{\R}$ &--& $\R\cup\{\infty\}$.
            	\end{tabular}
            \end{center}
                
            \vspace{1mm}
            
            In Section~\ref{Sec:ProofMainResults}, we will require that $p\neq 2,3,11$.

       \vspace{3mm}
        
	        \begin{center} \textbf{Acknowledgements} \end{center}

		This work started during the first and second named authors' visit to the Max Planck Institute for Mathematics in the Sciences. The authors thank the aforementioned institute for the generous hospitality and welcoming environment. The authors are grateful to Bernd Sturmfels for illuminating discussions. They also would like to thank the anonymous referee for several useful suggestions which improved the readability of this paper. The first author was supported by UK Research and Innovation under the Future Leaders Fellowship programme MR/S034463/2. The third author is supported by grant GYN-D9843-G0B1721N of the Fund for Scientific Research–Flanders.

\medskip

  {\it{ \bf UKRI Data Access Statement}}: the code associated to this article is publicly available at the link \eqref{eq:link_for_code}. Apart from this, there are no further associated research materials.

    \section{Background}\label{Sec:Background}
    
    In this section, we review the necessary background and preliminaries on invariant theory for binary forms and $(4,1)$-forms. We also recall some notions on Picard curves.
    
    \subsection{Invariant theory of binary forms}
    \label{Subsec:InvThy}
    
    In recall some facts and results from invariant theory of binary forms. For detailed treatments of invariant theory, we refer the reader to \cite{sturmfels2008algorithms,derksen2015computational,Draisma}. Fix an algebraically closed field $K$ such that $\mathrm{char}(K)=0$.

    \subsubsection{A group action on binary forms} 
    	 Let $n\geq{0}$ be a fixed integer. Let $A=K[a_{0},\dots,a_{n}]$ be the polynomial ring in $n+1$ variables. We view $A$ as a graded ring with the standard grading $\mathrm{deg}(a_{i})=1$.  Let $V_n$ be the $A$-submodule of $A[x,z]$ consisting of homogeneous polynomials in $x$ and $z$ of total degree $n$. The reductive group $ G \coloneqq \SL_2(K)$ acts (as a right action) on the $A$-module $V_n$ as follows: 
	\begin{equation} \label{eq:changeOfVariable}
	    g^{\sigma}(x,z) \coloneqq g(ax + bz, cx + dz), \quad \text{ for } g \in V_n \text{ and } \sigma = \begin{bmatrix} a & b \\ c & d	 \end{bmatrix} \in G.
	\end{equation}
	
	\begin{definition}{\bf{[Universal binary form]}}
	We define the \textit{universal} binary form $f$ of degree $n$ over $K$ as the binary form given by
	\begin{equation*}
		f(x,z) = a_{0}x^{n} + a_{1}x^{n-1}z+ \cdots +a_{n}z^{n}\in V_n.
	\end{equation*}	
    A binary form over $K$ of degree $n$ is obtained by specializing the coefficients to $K$.
	\end{definition}
	
	We obtain an action of $G$ on $A$ by sending $a_{i}$ to the coefficient $c_i(f^{\sigma})$ of $x^{n-i} z^{i}$ in $f^{\sigma}$ for $0 \leq i \leq n$. This means that $G$ acts on $A$ as follows:
	\begin{equation} \label{eq:Gaction}
		F^{\sigma}(a_0, \dots, a_n) \coloneqq F(c_0( f^{\sigma} ), \dots, c_n(f^{\sigma})), \quad \text{ for } F \in A \text{ and } \sigma \in G.
	\end{equation}
	
	\begin{example} When $n = 2$, from the definition
		\[
			f^{\sigma}(x,z)  = f(ax + bz, cx + dz), \quad \text{ for } \sigma = \begin{bmatrix} a & b \\ c & d	 \end{bmatrix} \in G
		\]
		we get
		\[
		 	 \begin{bmatrix} c_0(f^{\sigma}) \\	c_1(f^{\sigma})  \\c_2(f^{\sigma})  \end{bmatrix} =	 \begin{bmatrix} a^2 & ac & c^2\\ 2 ab  & ad + bc & 2cd \\ b^2 &  bd & d^2 \end{bmatrix}  \begin{bmatrix} a_0 \\	a_1  \\ a_2	 \end{bmatrix}.
		\]
		We can then see how $\sigma$ acts on the generators $a_0,a_1,a_2$ of $A = K[a_0, a_1, a_2]$ and we have
		\[
		a_0^{\sigma} = a^2 a_0 + ac a_1 + c^2 a_2, \quad a_1^{\sigma} = 2ab a_0 + (ad + bc) a_1 + 2cd a_2 \quad \text{ and } a_2^{\sigma} = b^2 a_0 + bd a_1 + d^2 a_2.
		\]

	\end{example}
 
    \begin{definition}{\bf{[Separable binary forms]}} A binary form is said to be \textit{separable} if its discriminant does not vanish.
    \end{definition}
 
	\begin{definition}{\bf{[Invariants of binary forms]}}
	    A homogeneous polynomial $F \in{A}$ is called \textit{$G$-invariant} if
		\[
	        F^{\sigma} = F \quad \text{ for all } \sigma \in G = \SL_2(K).
	    \]
	   We denote the graded ring generated by all homogeneous $G$-invariant polynomials by $A^{G}$. 
	\end{definition}

    \medskip
 
	\begin{remark}\label{rem:GL2vsSL2}
	Notice that there is also an action of $\GL_2(K)$ on $A$ as in \eqref{eq:Gaction}. Since $K$ is algebraically closed, a homogeneous polynomial $F \in A$ is $\SL_2$-invariant if and only if for any $\sigma \in \GL_2$ we have
	\[
	    F^\sigma = \det(\sigma)^{\deg(F)} F.
	\]
	We say that a homogeneous polynomial $F\in{A}$ is \textit{$\GL_{2}$-invariant} if it satisfies the identity above. This immediately implies that the ring of invariants for $\GL_{2}$ and $\SL_{2}$ are the same. We will use these interchangeably.
	
	\end{remark}
	
	\medskip
 
    Since $G=\SL_{2}(K)$ is reductive, it is a well known fact from invariant theory that the ring $A^G$ is finitely generated over $K$; see, for example, \cite[Corollary~2.2.11]{derksen2015computational}. To find generators of $A^G$ we can use the notion of transvectants, or {\it{\"{U}berschiebung}}, which we explain briefly here. Let $g \in{V_{m}}$ and $h\in{V_{n}}$ with $m \geq n$ and let $r$ be an integer with $0\leq r \leq{n}$. We define the bilinear map
	\begin{equation*}
	\angles{\cdot , \cdot}_{r}	: V_m \times V_n  \xrightarrow[]{} V_{m + n - 2r},  \quad (g,h)\mapsto \sum_{i=0}^{r} (-1)^i\binom{r}{i} \dfrac{\partial^{r}{g}}{\partial^{r-i}{x} \ \partial^{i}{z}} \ \dfrac{\partial^r{h}}{\partial^{i}{x} \ \partial^{r-i}{z}} .
	\end{equation*}
	It turns out that this map is $G$-invariant, i.e.,
	\[
		\angles{g^{\sigma}, h^{\sigma}}_{r} = \angles{g,h}_{r}, \quad \text{ for } (g,h) \in V_m \times V_n \text{ and } \sigma \in G.
	\]
	The quantity $\angles{g,h}_{r}$ is called the $r$-th \textit{transvectant} of $g$ and $h$. Notice that when $n = m = r$, the transvectant $\angles{g,h}_{r}$ has degree $0$ in $x,z $ so  $\angles{g,h}_{r} \in A^G$. A theorem by Gordan \cite{Gordan1868} shows that one can obtain all invariants by taking (iterated if necessary) transvectants of $f$, so this makes the set generators of $A^G$ explicit, at least in theory. In practice, the generators have only been found for binary forms of degrees up to $10$; see \cite[Chapter~1]{Draisma} or \cite[Example~2.1.2]{derksen2015computational}.

	\begin{example}
		
	\begin{enumerate} \label{ex:InvsOfSmalln}
		\item For binary quadrics the invariant ring $A^G$ is generated by 
		$\angles{f,f}_{2}$, which is equal to the discriminant of $f$.  
		
		\item For binary cubics, the invariant ring is again generated by the discriminant
		\[
		    \Delta = \angles{\angles{f, f}_2,\angles{f, f}_2}_2.
		\]
		
		\item For binary quartics $n=4$, the invariant ring $A^G$ is generated by two algebraically independent invariants $I_{2}$ and $I_{3}$ of degrees $2$ and $3$ respectively. In terms of transvectants, they are given by
		\[
		    I_2 = \angles{f,f}_4 \quad \text{and} \quad I_3 = \angles{f,\angles{f,f}_2}_4.
		\]
		
		\item For binary quintics, the invariant ring $A^G$ is generated by four generators $I_{4}$, $I_{8}$, $I_{12}$ and $I_{18}$ of degrees $4$, $8$, $12$ and $18$, respectively. These are not algebraically independent, as $I_{18}^2$ can be expressed in terms of the other three. We refer the reader to \cref{subsection:quinticInvariants} for a more detailed exposition in this case.
	\end{enumerate}

	\end{example}

    \subsubsection{Invariants of binary quintics} \label{subsection:quinticInvariants}
				    
		    As discussed in Example~\ref{ex:InvsOfSmalln}, the ring of invariants for quintic binary forms is generated by four elements $I_{4}$, $I_{8}$, $I_{12}$ and $I_{18}$ of degrees $4$, $8$, $12$ and $18$, respectively. Let $f$ be the universal quintic
		    \[
		    f(x,z) = a_0 x^5 + a_1 x^4 z + a_2 x^3 z^2 + a_3 x^2 z^3 + a_4 x z^4 + a_5 z^5.
		    \]
		    The four invariants $I_4$, $I_8$, $I_{12}$ and $I_{18}$ can be obtained by taking transvectants of $f$ in the following way:
		    \begin{equation}\label{eq:I_invariants}
		    \begin{aligned}
				I_{4}&=\angles{f^2,f^2}_{10},  & I_{8}&=\angles{f^4,f^4}_{20},\\
				I_{12}&=\angles{f^6,f^6}_{30}, & I_{18}&=\angles{\angles{f^5,f^6}_{10},f^7}_{35}. 
			\end{aligned}  
		    \end{equation}
		    Note that these transvectants are different from the ones considered in \cite[Section~4.4]{Draisma}\footnote{The universal binary form used in the literature is of the shape $f = \sum_{j = 0}^{n} \binom{n}{j} a_j x^{n-j} z^{j}$. In this paper we use $f = \sum_{j = 0}^{n} a_j x^{n-j} z^{j}$, since in our field $K$ some binomial coefficients $\binom{n}{j}$ might have a positive valuation. Since $\mathrm{char}(K) = 0$, the invariant theory remains unchanged.}. They still yield generators however, as one easily checks using the Poincar\'{e} series. The discriminant $\Delta$ of $f$ can be expressed in terms of the previous invariants as
		    \begin{equation} \label{eq:Discrim}
				    \Delta= c_{0}I_{4}^2 + c_{1} I_{8}
		    \end{equation}
		    where $c_0, c_1$ are the constants
				\[
				    c_{0} =-1/2746158938062848000000, \quad  c_{1} = 1/46987474647852089270599680000000.
				\]
			Note that $\deg(\Delta) = 8$.
			
			In \cref{Sec:tropicalInvariants}, we will assume that $K$ is a complete non-archimedean field and associate a metric tree to a binary form over $K$. This tree is invariant under the action of $\GL_{2}$. For quintics, there are three different types, see Figure~\ref{fig:tropicalTypes}. The invariants $I_4,I_8,I_{12},I_{18}$ and the discriminant $\Delta$ are however not enough to distinguish between the different types, as we will see in 
			\cref{ex:counterExampleTreeType}. In order to distinguish tree types, we need to introduce another invariant, which we call the \emph{$H$-invariant}.

			\begin{definition}{\bf{[$H$-invariant]}}
				The $H$-\emph{invariant} of the quintic $f = a_0 x^5 + a_1 x^4 z + \dots + a_5 z^5$ is defined as follows:
				\begin{equation} \label{eq:H_inv}
						H =  \beta I_{12} - 396  \alpha^3 I_{4}^{3},
				\end{equation}
				where $\alpha$ and $\beta$ are given by
				\begin{align*}
					\alpha &= 2^{-17} \cdot 3^{-7} \cdot 5^{-3} \cdot 7^{-1},\\
					\beta &= 2^{-50} \cdot 3^{-27} \cdot 5^{-14} \cdot 7^{-7} \cdot 11^{-4} \cdot 13^{-3} \cdot 17^{-1} \cdot 19^{-1} \cdot 23^{-1} \cdot 29^{-1}.
				\end{align*}
			\end{definition}
			
			This invariant is of degree $12$ and constructed to satisfy the following property:
			\[
			       H(0,0,1,0,a_4,a_5) \quad \text{ is equal to the discriminant of the cubic } \quad x^{3} + a_4x + a_5.
			\]
			The motivation behind creating this invariant is as follows. For binary quintics, there are three unmarked tree types, see Figure~\ref{fig:tropicalTypes}. For the two non-trivial ones, we  choose a non-trivial vertex of valency two. By applying a projective linear transformation, we can ensure that the roots are grouped as $\{\infty,\lambda_{1}\}$ and $\{0,1,\lambda_{2}\}$, where $v(\lambda_{1})<v(\lambda_{2})=0$. After scaling the binary form so that the smallest valuation of its coefficients is zero, the coefficients of $x^5$ and $x^4z$ have positive valuation. 
            This binary form hence reduces to a cubic and the discriminant of this cubic distinguishes between the two remaining types.

    \subsubsection{Invariants of \texorpdfstring{$(4,1)$}~-forms}\label{subsection:InvariantsFourOneForms}
    
    Let $B = K[b_0, \dots, b_4, c_0, c_1]$ be a polynomial ring. For an integer $n \geq 0$, we denote by $W_n$ the space of homogeneous polynomials in $B[x,z]$ of degree $n$. Here, we are interested in $(4,1)$-binary forms, which are the elements $(q, \ell)$ of $W_4 \oplus W_1$. There is a natural group action of $G = \SL_2(K)$ on $(4,1)$-forms as follows
		\[
			(q(x,z), \ell(x,z))^{\sigma} = ( q^{\sigma}(x,z), \ell^{\sigma} (x,z)) = ( q(\sigma(x,z)), \ell(\sigma (x,z))).
		\]
		
	\begin{definition}{\bf{[Universal $(4,1)$-form]}}
	We define the \textit{universal} $(4,1)$-form $(q,\ell)$ as
	\[
	q(x,z) = b_0x^4 + b_1 x^3 z + \dots + b_4 z^4 \quad \text{and} \quad \ell(x,z) = c_0 x + c_1 z.
	\]
	\end{definition}	
		
    The action of $G$ on $(q,\ell)$ defines an action of $G$ on $B$ where $\sigma \in G$ acts on $B$ by sending $ b_0, \dots, b_4, c_0, c_1$ to the coefficients of $(q,\ell)^{\sigma}$. We denote by $B^G$ the ring of $G$-invariant polynomials in $B$. This ring is again finitely generated since $G$ is reductive. It has $5$ generators $j_2, j_3, j_5,j_6$ and $j_9$ of respective degrees $2,3,5,6$ and $9$, and they can be obtained using transvectants as follows:
    \begin{equation}\label{eq:j_Inv}
    \begin{aligned}
        j_{2}&=\angles{q,q}_{4},		    & j_{3}&=\angles{\angles{q,q}_{2},q}_{4}, \\  
        j_{5}&=\angles{q,\ell^4}_{4},       & j_{6}&=\angles{\angles{q,q}_{2},\ell^{4}}_{4},\\
        j_{9}&=\angles{\angles{q,\angles{q,q}_{2}}_1, \ell^6}_6. 
    \end{aligned}        
    \end{equation}
    We refer the reader to \cite[Section~5.4.1]{Draisma} for more details.
    
     \vspace{0.3cm}
      
    The invariants $I_4,  I_8, I_{12}, I_{18}, \Delta, H$ for binary quintics and the invariants $j_2, j_3, j_5, j_6, j_9$ for $(4,1)$-forms are polynomials with rational coefficients.
    Scaling these invariants, we may assume that their coefficients are integer and are coprime; recall that $\mathrm{char}(K) = 0$.
    The reason behind this is that we are working with a valued field $K$ in which integer scalars might have a positive valuation. This justifies the following assumption:
      
    \begin{assumption}
    We scale all the invariants so that the greatest common divisor of their coefficients is $1$, and we do not change the notation. These scaled invariants are computed in (\ref{eq:link_for_code}), and are the ones used in our results.
	\end{assumption}
    
    Write $B_{4,1}=B$ for the coefficient ring corresponding to $(4,1)$-forms and $B_{5}$ for the coefficient ring corresponding to quintics. The universal $(4,1)$-form $(q,\ell)$ gives a canonical quintic $f=q\cdot\ell$. From this, we obtain an injective ring homomorphism 
    \begin{equation*}
        B_{5}\to{B_{4,1}}.
    \end{equation*}
    This ring homomorphism is $G = \SL_{2}(K)$-equivariant, so we obtain an induced injective ring homomorphism
    \begin{equation*}
        B_{5}^{G}\to{B^{G}_{4,1}}
    \end{equation*}
    of invariants. We identify $B^{G}_{5}$ with its image, so that we have an inclusion $B_{5}^{G}\subset{}{B^{G}_{4,1}}$. The generators from Section~\ref{subsection:quinticInvariants} can then be expressed in terms of the $j_{i}$ explicitly as follows:
        \begin{equation}
        \begin{aligned}\label{eq:InvariantsInTermsofJs}
          I_4     &=  \frac{2}{3}j_2 j_6  - \frac{1}{2} j_3 j_5, \\
          I_8     &= \frac{14}{9}j_2^2 j_6^2 + \frac{22}{27} j_2^3 j_5^2 + \frac{5}{27} j_3^2 j_5^2 - \frac{14}{9} j_2 j_3 j_5 j_6, \\
          I_{12}  &= \frac{4400}{243}  j_2^3  j_6^3  - \frac{11}{243} j_3^2  j_6^3  - \frac{242}{9}  j_2^2  j_3  j_5  j_6^2  + \frac{2479}{81}  j_2^4  j_5^2  j_6  + \frac{692}{81} j_2 j_3^2 j_5^2 j_6\\
                   & \ \ \ + \frac{7156}{243} j_2^3 j_3 j_5^3  - \frac{92}{243} j_3^3 j_5^3,  \\ 
          I_{18}  &= - \frac{625}{729} j_2^6  j_5^3  j_9  - \frac{512}{729} j_2^3 j_3^2  j_5^3  j_9   - \frac{4}{729} j_2^3 j_3 j_6^3 j_9 + \frac{1}{729} j_3^3 j_6^3 j_9 + 
          \frac{1}{3} j_2^4 j_3 j_5^2 j_6 j_9 \\ 
          & \ \ \ + \frac{4}{243} j_2^5 j_5 j_6^2  j_9 - \frac{1}{243} j_2^2  j_3^2  j_5  j_6^2  j_9;
        \end{aligned}
        \end{equation}
      see the computations in \eqref{eq:link_for_code}. For $\Delta$ and $H$, we use \eqref{eq:Discrim} and \eqref{eq:H_inv}. 

    Finally, we summarize the invariants used in Theorems \ref{thm:tropicalTypes}, \ref{thm:reductionTypes} and \ref{thm:edgeLengths}, and their respective degrees as follows:

  	  \begin{table}[H]
            \centering
                \begin{tabular}{|c|c|c|c|c|c|c|c|}
                \hline
                    \textbf{Invariant} & $I_4$ & $I_8$ & $I_{12}$ & $I_{18}$ & $\Delta$ & $H$  \\
                \hline     
                    \textbf{Degree}    & $4$ & $8$ & $12$ & $18$ & $8$ & $12$\\
                \hline
                \end{tabular} \hspace{5mm}
                \begin{tabular}{|c|c|c|c|c|c|}
                \hline
                    \textbf{Invariant} & $j_2$ & $j_3$ & $j_{5}$ & $j_{6}$ & $j_{9}$ \\
                \hline     
                    \textbf{Degree}    & $2$   & $3$   &   $5$   & $6$ & $9$ \\
                \hline
                \end{tabular}
            \caption{Invariants of binary quintics (on the left) and (4,1)-forms (on the right) together with their degrees.}
        \label{tab:my_label}
    \end{table}

    \subsection{Picard curves}
    
    Let $K$ denote a field of characteristic $0$. A Picard curve $X$ over $K$ is a smooth projective curve of genus $3$ in $\mathbb{P}^{2}$ given by an equation of the form 
      \[
        y^3\ell(x,z) = q(x,z),
    \]
    where $q$ and $\ell$ are homogeneous polynomials of degrees $4$ and $1$ respectively in $K[x,z]$. Here the smoothness of $X$ is equivalent to the separability of the quintic $q(x,z)\cdot\ell(x,z)$. For a Picard curve, the five distinct roots of this quintic are the branch points of the $\Z/3\Z$-covering of $\mathbb{P}^{1}$ induced by the rational map
    \[
                \pi \colon X \DashedArrow \P^1, \quad [x:y:z] \mapsto [x:z].
    \]
    Either using the fact that $X$ is smooth or an explicit computation, we find that it gives a morphism, so that it is defined everywhere. 
    The Galois group of $\pi$ is $G=\mathbb{Z}/3\mathbb{Z}$ and it acts as follows on $X$. Let $\zeta\in{K}$ be a primitive third root of unity. There is then an automorphism 
      \[
      \alpha:[x:y:z]\mapsto [x:\zeta{y}:z]
      \]
      of order three and the quotient map $X\to X/\langle \alpha \rangle =\mathbb{P}^{1}$ is $\pi$. For each of the five points $P_{i}$ in the zero set of the quintic $q(x,z)\cdot\ell(x,z)$, there is a unique point $Q_{i}$ lying over $P_{i}$. The four points $Q_{1}$,\dots,$Q_{4}$ corresponding to the quartic $q(x,z)$ differ from the remaining point $Q_{5}$ in the sense that $\alpha$ acts by $\zeta$ on the tangent spaces of $Q_{1}$,\dots,$Q_{4}$ and as $\zeta^{2}$ on the tangent space of $Q_{5}$. This gives a second, and equivalent, definition of a Picard curve as a $\Z/ 3\Z$-covering of $\mathbb{P}^{1}$ branched over five points with a specified action of $\mathbb{Z}/3\Z$ on the tangent spaces of the ramification points; see \cite{AchterPries}.  This rigidification can also be given in terms of differential forms, see \cite[Section~5]{Clery_vanderGeer_21}.

    \section{Tropical invariants for general group actions}
    \label{Sec:tropicalInvariants}
    
    In this section, we define the notion of a set of tropical invariants for general group actions. The main underlying idea is that they separate orbits that have the same limit points in some toric compactification. For binary forms, this gives us the notion of a set of tropical invariants. For Picard curves, there are two interesting types of binary forms: quintics and $(4,1)$-forms. We will see in \cref{Sec:ProofMainResults} that we can give a tropical set of invariants for these forms.

    We assume throughout this section that $K$ is a complete, non-archimedean and algebraically closed field $K$ of characteristic $0$ with valuation ring $R$ and residue field $k$.
    
    \subsection{Compactifications and tropicalizations}\label{sec:CompactificationsTropicalizations}
    
    Let $U$ be an irreducible variety over $K$ with Berkovich analytification $U^{\mathrm{an}}$. We can define tropical separating sets for partitions of subsets of $U^{\mathrm{an}}$ as follows.

    \begin{definition}{\bf{[Separating set]}}
    Let $U^{\mathrm{an}}\supset{S}=\bigsqcup S_{i}$ be a partition of a subset $S$ of $U^{\mathrm{an}}$ into disjoint subsets and let $\phi: U\to{X(\Delta)}$ be an embedding of $U$ into a toric variety $X(\Delta)$ with fan $\Delta$. We say that $\phi$ is \textit{separating} for the partition if $\mathrm{trop}(\phi(S_{i}))\cap\mathrm{trop}(\phi(S_{j}))=\emptyset$ for $i\neq{j}$, where $\trop(\cdot)$ is the natural tropicalization map associated to the toric variety $X(\Delta)$. If the embedding is given by the global sections of a line bundle on $U$ and $\phi$ is separating for a partition, then we call the sections a \textit{tropical separating set} for the partition.  If the sections of the line bundle give a morphism to $X(\Delta)=\mathbb{P}(a_{1},\dots,a_{n})$, then we say that they are a \textit{projective tropical separating set}. 
    \end{definition}

    \begin{example}
        If we take $S=\{P_{1}\}\cup\{P_{2}\}$ for two points $P_{1},P_{2}\in U^{\mathrm{an}}$, then the proof of \cite[Theorem~1.1]{payne2009analytification} shows that we can find an embedding that separates these points.
    \end{example}

    The partitions of $U^{\mathrm{an}}$ we are interested in are given using compactifications of $U$ as follows. Let $\mathcal{X}\to{\mathrm{Spec}(R)}$ be a proper flat $R$-scheme and let $U\to{\mathcal{X}}$ be an open immersion. Consider the reduction map 
    \begin{equation}\label{eq:ReductionMap}
        \mathrm{red}:U^{\mathrm{an}}\to\mathcal{X}_{s}
    \end{equation}
    as in \cite[Section~4]{GublerGuideToTrop}, where $\mathcal{X}_{s}$ is the special fiber of $\mathcal{X}$. We now obtain a partition of $U^{\mathrm{an}}$ by taking the inverse image under $\mathrm{red}(\cdot{})$ of a partition of $\mathcal{X}_{s}$.

\subsection{Tropical invariants}

    We start with the definition of tropical weighted projective space.

    \begin{definition}{\bf{[Tropical weighted projective space]}} \label{def:weightedTropiclaProjSpace}
    Let $(\T\A^n)^* = \wb{\R}^n\setminus\{(\infty,\dots,\infty)\}$, the punctured tropical affine plane. For a fixed $n$-tuple of natural numbers $(a_1,\dots,a_n)$, we define an action of $\R$ on $(\T\A^n)^*$ as follows:
    \[
    \lambda\odot(x_1,\dots,x_n)\coloneqq (x_1,\dots,x_n)+\lambda\cdot(a_{1},\dots,a_{n}).
    \]
    The (set-theoretic) quotient of $(\T\A^n)^*$ by this action is called the \textit{tropical weighted projective space}, and is denoted by $\T\P(a_1,\dots,a_n)$.
    \end{definition}
    
    Now, let $U/K$ be an irreducible variety and let $G/K$ be a group scheme acting on $U$. Suppose that there exists a geometric quotient $V=U/G$ for $U$ and $G$. The notion of a tropical separating from Section~\ref{sec:CompactificationsTropicalizations} now allows us to define a set of tropical invariants.  
    
    \begin{definition}
    {\bf{[Tropical invariants]}} Let $V^{\mathrm{an}}\supset{S}=\bigsqcup {S_{i}}$ be a partition of a subset $S$ of $V^{\mathrm{an}}$. A set of \textit{tropical invariants} for the partition is a tropical separating set for the partition. If the tropical invariants give an embedding into a weighted projective space $\mathbb{P}(a_{1},\dots,a_{n})$, then we call this a set of \textit{tropical projective invariants}. The image of the tropicalization is the space $\T\P(a_1,\dots,a_n)$.
    \end{definition}

    Since points of $V^{\mathrm{an}}$ correspond to orbits in $U^{\mathrm{an}}$, we find that a set of tropical invariants separates a partition of a subset of $U^{\mathrm{an}}$ that is stable under the action of $G$.

    \begin{example}
        Consider the ring $K[x,y]$. The group $S_{2}\simeq{}\mathbb{Z}/2\mathbb{Z}$ acts on $K[x,y]$ by $x\mapsto{y}$ and $y\to{x}$. The invariant ring is $K[xy,x+y]$ and, on the level of schemes, the map $\mathrm{Spec}(K[x,y])\to\mathrm{Spec}(K[xy,x+y])$ gives the geometric quotient. Consider the standard embedding of $\mathbb{A}^{2}=\mathrm{Spec}(K[x,y])$ into $\mathbb{P}^{2}$. The group action then extends to $\mathbb{P}^{2}$. The boundary is isomorphic to $\mathbb{P}^{1}$ and we decompose $\mathbb{P}^{2}$ into 
        \[
            \overline{S}_{0} = \{[x:y:z] \mid z\neq{0}\},\ \overline{S}_{1}=\{[1:1:0]\} \text{ and } \overline{S}_{2}=\{[x:y:0] \mid x\neq{y}\}.
        \] 
        This induces a partition
        \[
            S_{0}\sqcup{}S_{1}\sqcup{S_{2}}=\Spec(K[xy,x+y])^{\mathrm{an}}
        \] 
        since the group action preserves the $\overline{S_{i}}$. We now note that the standard embedding
        \[\mathrm{Spec}(K[xy,x+y])=\mathbb{A}^{2}\to\mathbb{P}^{2}\]
        is not tropically separating since $\mathrm{char}(K) \neq 2$. If we however consider the embedding into $\mathbb{P}^{3}$ given by $\{xy,x+y,(x-y)^2,1\}$, then this does form a tropical set of invariants. 
    \end{example}

      We will see many more examples of tropical invariants in Sections \ref{sec:ModuliStable} and \ref{sec:TreesTropicalBinaryForms}. 
      
 \subsection{Moduli of tropical curves}
 \label{sec:ModuliStable}
 
 We now review moduli of stable curves of genus zero. The functor $\mathcal{M}_{0,n}: \bf{Sch} \to  \bf{Sets}$ from the category of schemes to the category of sets, defined by 
 \begin{equation}\label{eq:FunctorDefinition}
     \mathcal{M}_{0,n}(S)=\{\text{smooth curves } C\to{S} \text{ of genus zero with }n \text{ marked points}\}/\sim{}
 \end{equation}
 for $n\geq{3}$ is representable by a scheme $M_{0,n}$. Note that there are indeed no non-trivial elements of $\mathrm{Aut}(\mathbb{P}^{1})=\mathrm{PGL}_{2}$ that fix three or more elements of $\mathbb{P}^{1}$. We refer to $M_{0,n}$ as the moduli space of smooth $n$-marked curves of genus zero. There is a natural compactification of this space, given by replacing smooth curves by stable curves in \eqref{eq:FunctorDefinition}. We denote the corresponding scheme by $\overline{M}_{0,n}$. It is proper and flat over $\mathrm{Spec}(\Z)$. The boundary locus $\Delta=\overline{{M}}_{0,n}\backslash{M_{0,n}}$ corresponds to non-smooth stable curves.
 
 The group $S_{n}$ acts on the above moduli spaces by permuting the marked points. Since $S_{n}$ is finite and $M_{0,n}$ and $\overline{M}_{0,n}$ are quasi-projective, there exist geometric quotients $M_{0,n}/S_{n}$ and $\overline{M}_{0,n}/S_{n}$ which fit into a commutative diagram 
 \begin{equation*}
     \begin{tikzcd}
{M_{0,n}} \arrow[r] \arrow[d] & {\overline{M}_{0,n}} \arrow[d] \\
{M_{0,n}/S_{n}} \arrow[r]           & {\overline{M}_{0,n}/S_{n}}          
\end{tikzcd}
 \end{equation*} The scheme ${\overline{M}_{0,n}/S_{n}}$ is again proper over $\mathbb{Z}$. The geometric points of ${M_{0,n}/S_{n}}$ (resp. ${\overline{M}_{0,n}/S_{n}}$) can be identified with smooth (resp. stable) curves of genus $0$ with $n$ unordered points.
 
 We now describe their abstract tropicalizations as in \cite[Section~4]{DanCaporasoPayne}. We start with the original set-up without quotients. Consider phylogenetic trees with $n$ marked leaves, together with a length function $\ell$ on the non-leaves. These are the points of the tropical moduli space $M_{0,n}^{\mathrm{trop}}$. It can be given the structure of a generalized cone complex as in \cite[Section~2]{DanCaporasoPayne}. It can also be given as the tropicalization of the Pl\"{u}cker map as in \cite[Theorem 6.4.12]{MaclagenSturmfels}. We obtain an abstract tropicalization map 
 \begin{equation*}
     \mathrm{trop}:M^{\mathrm{an}}_{0,n}\to{M^{\mathrm{trop}}_{0,n}}
 \end{equation*}
by the following procedure. A point  $P\in{M^{\mathrm{an}}_{0,n}}$ can be represented by an $L$-valued point of $M_{0,n}$, and we use the natural map
\begin{equation}\label{eq:ValuativeLimit}
    M_{0,n}(L)\to{\overline{M}_{0,n}(L)}=\overline{M}_{0,n}(R_{L})
\end{equation}
to obtain a stable model $\mathcal{X}\to\mathrm{Spec}(R_{L})$ with special fiber $\mathcal{X}_{s}\to\mathrm{Spec}(k_{L})$. Here, $L$ is a complete valued field extension of $K$ with valuation ring $R_{L}$ and residue field $k_{L}$. Let $T$ be the marked dual intersection graph of $\mathcal{X}_{s}$, together with the natural edge length function $\ell$ induced from  $\mathcal{X}$. Note that there is a unique leaf for every marked point. We set $\mathrm{trop}(P)=(T,\ell)$.

We can also express the tropicalization of a point in $M_{0,n}^{\an}$ in terms of Berkovich skeleta as in \cite{Baker_Payne_Rabinoff_13}. Namely, a point $P\in M_{0,n}^{\an}$ corresponds to a marked curve $(\mathbb{P}^{1},(P_{1},...,P_{n}))$, where $P_{i}\in\mathbb{P}^{1}(L)$ are distinct $L$-valued points. The minimal skeleton $\Sigma$ of this marked curve as in \cite[Section 2.6]{Baker_Payne_Rabinoff_13}\footnote{In \cite{Baker_Payne_Rabinoff_13}, these are referred to as punctured curves.} is then a marked metric tree whose leaves correspond to the points $P_{i}$.
The underlying graph of this marked metric tree is the same as $T$ and it gives rise to these same edge length function $\ell$. By abuse of notation, we will write $\Sigma=(T,\ell)$, where we either view $(T,\ell)$ as a marked finite graph with a length function, or a marked metric tree.

\begin{definition}\label{def:MarkedMetricTree}
    Let $P\in M_{0,n}^{\an}$ for $n\geq{3}$ with induced marked curve $(\mathbb{P}^{1},(P_{1},...,P_{n}))$. We call $\Sigma=(T,\ell)$ the marked metric tree associated to $P$. 
\end{definition}

\begin{remark}
We note that the edge lengths of the marked metric tree $\Sigma$ are independent of the chosen representative of the isomorphism class corresponding to $P$. That is, if we have an isomorphism of marked curves $(\mathbb{P}^{1},(P_{1},...,P_{n}))\to (\mathbb{P}^{1},(\sigma(P_{1}),...,\sigma(P_{n})))$, then this induces a morphism of marked metric trees $\Sigma\to \sigma(\Sigma)$ that is an isometry outside the leaves. This, for instance, follows from the results in \cite[Section 2]{Baker_Payne_Rabinoff_13}. 
\end{remark}

Let $\mathcal{P}([n])$ be the power set of $[n]=\{1,...,n\}$ and let $\Xi\subset\mathcal{P}([n])$ be a partition. We write $S_{\Xi}\subset S_{n}$ for the subgroup of permutations that preserve $\Xi$. Note that $S_{\Xi}$ is isomorphic to a product of symmetry groups $S_{k}$ for $k\leq{n}$. We now define an action of $S_{\Xi}$ on $M^{\mathrm{trop}}_{0,n}$. We view an element of $M^{\mathrm{trop}}_{0,n}$ as a metric tree $\Sigma=(T,\ell)$, together with an injection $i:\{1,\dots,n\}\to{L(T)}$, where $L(T)$ is the set of (infinite) leaves of $T$. By permuting the leaves indexed by the partition $\Xi$, we obtain an action of $S_{\Xi}$ on $M^{\mathrm{trop}}_{0,n}$, where the latter is viewed as an object in the category of generalized cone complexes. We now note that categorical quotients in the category of generalized cone complex exist since arbitrary finite colimits exist, see \cite[Remark~2.6.1]{DanCaporasoPayne}. From this, we obtain the following definition of $M^{\mathrm{trop}}_{0,n}/S_{\Xi}$.

 \begin{definition}\label{def:TropicalModuliSpace}
     Let $\Xi\subset\mathcal{P}([n])$ be a partition of $[n]=\{1,...,n\}$ and consider the induced action of $S_{\Xi}$ on $M^{\mathrm{trop}}_{0,n}$. The space $M^{\mathrm{trop}}_{0,n}/S_{\Xi}$ is  the categorical quotient of $M^{\mathrm{trop}}_{0,n}$ under the action of $S_{\Xi}$. A $\Xi$-tree-type is an equivalence class of a tree type in $M^{\mathrm{trop}}_{0,n}$ under the action of $S_{\Xi}$. 
 \end{definition}

\begin{example}
    Suppose that $\Xi=\{\{1,...,n\}\}$. The corresponding moduli space $M^{\mathrm{trop}}_{0,n}/S_{n}$ parametrizes tropical trees with $n$ unmarked points. 
\end{example}

 \begin{example}
 There are three non-trivial marked tree types for $n=4$, giving three cones $\mathbb{R}_{\geq{0}}$. The types with edge lengths zero are all identified, so we glue these cones together to obtain $M^{\mathrm{trop}}_{0,4}$, which is a standard tropical line. The three cones are permuted by $S_{4}$, giving the quotient $M^{\mathrm{trop}}_{0,4}/S_{4}=\mathbb{R}_{\geq{0}}$.
 \end{example}
 
 \begin{example}\label{exa:TreesFiveLeaves} Let $n=5$. Then, there are exactly three unmarked types: $\mathrm{I}$, $\mathrm{II}$ and $\mathrm{III}$; they are depicted in Figure~\ref{fig:tropicalTypes}.
 Type~$\mathrm{III}$ corresponds to a folded positive orthant of dimension $2$, since the automorphism group of the underlying graph is $\Z/2\Z$. Figure~\ref{fig:M_trop_05} represents the space $M^{\mathrm{trop}}_{0,5}/S_{5}$.

 \begin{figure}[H]
    \centering
    \scalebox{0.6}{ 
      \tikzset{every picture/.style={line width=0.75pt}}
        
        \begin{tikzpicture}[x=0.75pt,y=0.75pt,yscale=-1,xscale=1]
            \draw  [color={rgb, 255:red, 0; green, 0; blue, 0 }  ,draw opacity=0 ][fill={rgb, 255:red, 155; green, 155; blue, 155 }  ,fill opacity=0.6 ] (390,310) -- (170,310) -- (390,100) -- cycle ;
            \draw [line width=2.25]    (390,310) -- (170,310) ;
            \draw  [fill={rgb, 255:red, 0; green, 0; blue, 0 }  ,fill opacity=1 ] (165,310) .. controls (165,307.24) and (167.24,305) .. (170,305) .. controls (172.76,305) and (175,307.24) .. (175,310) .. controls (175,312.76) and (172.76,315) .. (170,315) .. controls (167.24,315) and (165,312.76) .. (165,310) -- cycle ;
            \draw    (138.12,315.6) -- (156.23,335) ;
            \draw    (156.23,296.21) -- (138.12,315.6) ;
            \draw    (138.12,315.6) -- (138.12,289.74) ;
            \draw    (120,296.21) -- (138.12,315.6) ;
            \draw    (138.12,315.6) -- (120,335) ;
            \draw  [fill={rgb, 255:red, 0; green, 0; blue, 0 }  ,fill opacity=1 ] (135.1,315.6) .. controls (135.1,313.82) and (136.45,312.37) .. (138.12,312.37) .. controls (139.78,312.37) and (141.14,313.82) .. (141.14,315.6) .. controls (141.14,317.39) and (139.78,318.84) .. (138.12,318.84) .. controls (136.45,318.84) and (135.1,317.39) .. (135.1,315.6) -- cycle ;
            \draw    (265.53,131.67) -- (315.12,131.67) ;
            \draw  [fill={rgb, 255:red, 0; green, 0; blue, 0 }  ,fill opacity=1 ] (263.05,131.67) .. controls (263.05,129.98) and (264.16,128.62) .. (265.53,128.62) .. controls (266.9,128.62) and (268.01,129.98) .. (268.01,131.67) .. controls (268.01,133.36) and (266.9,134.73) .. (265.53,134.73) .. controls (264.16,134.73) and (263.05,133.36) .. (263.05,131.67) -- cycle ;
            \draw  [fill={rgb, 255:red, 0; green, 0; blue, 0 }  ,fill opacity=1 ] (287.84,131.67) .. controls (287.84,129.98) and (288.95,128.62) .. (290.32,128.62) .. controls (291.69,128.62) and (292.8,129.98) .. (292.8,131.67) .. controls (292.8,133.36) and (291.69,134.73) .. (290.32,134.73) .. controls (288.95,134.73) and (287.84,133.36) .. (287.84,131.67) -- cycle ;
            \draw  [fill={rgb, 255:red, 0; green, 0; blue, 0 }  ,fill opacity=1 ] (312.64,131.67) .. controls (312.64,129.98) and (313.75,128.62) .. (315.12,128.62) .. controls (316.49,128.62) and (317.6,129.98) .. (317.6,131.67) .. controls (317.6,133.36) and (316.49,134.73) .. (315.12,134.73) .. controls (313.75,134.73) and (312.64,133.36) .. (312.64,131.67) -- cycle ;
            \draw    (290.32,107.23) -- (290.32,131.67) ;
            \draw    (250.65,113.34) -- (265.53,131.67) ;
            \draw    (265.53,131.67) -- (250.65,150) ; 
            \draw    (315.12,131.67) -- (330,150) ;
            \draw    (330,113.34) -- (315.12,131.67) ;
            \draw    (278.41,339.71) -- (309.09,339.71) ;
            \draw  [fill={rgb, 255:red, 0; green, 0; blue, 0 }  ,fill opacity=1 ] (275.34,339.71) .. controls (275.34,337.89) and (276.71,336.42) .. (278.41,336.42) .. controls (280.1,336.42) and (281.48,337.89) .. (281.48,339.71) .. controls (281.48,341.52) and (280.1,342.99) .. (278.41,342.99) .. controls (276.71,342.99) and (275.34,341.52) .. (275.34,339.71) -- cycle ;
            \draw  [fill={rgb, 255:red, 0; green, 0; blue, 0 }  ,fill opacity=1 ] (306.02,339.71) .. controls (306.02,337.89) and (307.4,336.42) .. (309.09,336.42) .. controls (310.79,336.42) and (312.16,337.89) .. (312.16,339.71) .. controls (312.16,341.52) and (310.79,342.99) .. (309.09,342.99) .. controls (307.4,342.99) and (306.02,341.52) .. (306.02,339.71) -- cycle ; 
            \draw    (260,320) -- (278.41,339.71) ;
            \draw    (278.41,339.71) -- (260,359.42) ;
            \draw    (309.09,339.71) -- (327.5,359.42) ;
            \draw    (327.5,320) -- (309.09,339.71) ;
            \draw    (309.09,339.71) -- (327.5,339.71) ;
        \end{tikzpicture}
         }
    \caption{The space $M^{\mathrm{trop}}_{0,5}/S_{5}$}
    \label{fig:M_trop_05}
\end{figure}
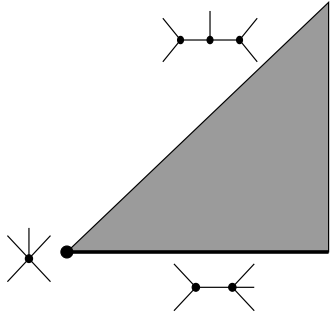

 \end{example}

 \begin{example}

     Suppose that $\Xi=\{\{1,2,...,n-1\},\{n\}\}$. In this case, we will write $S_{\Xi}=S_{n-1}$ for brevity. The corresponding tropical moduli space $M^{\mathrm{trop}}_{0,n}/S_{n-1}$ parametrizes tropical trees with one marked point.   

   For instance, consider the case where $n=5$ and $\Xi=\{\{1,2,...,4\},\{5\}\}$, giving rise to an action of $S_{4}$ on $M^{\mathrm{trop}}_{0,5}$. The quotient 
   $M^{\mathrm{trop}}_{0,5}/S_{4}$ is the moduli space of phylogenetic trees with $4$ unmarked leaves and $1$ marked leaf, see Figure~\ref{fig:M_trop_041}.
   There are five different corresponding phylogenetic tree types: $\mathrm{I}$, $\mathrm{II}.1$, $\mathrm{II}.2$, $\mathrm{III}.1$ and $\mathrm{III}.2$; see Figure~\ref{fig:fourOneTreeTypes}.

 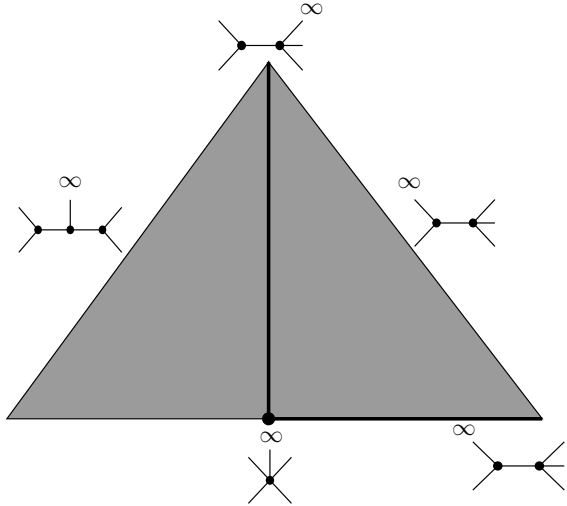
\begin{figure}[ht]
    \centering
    \scalebox{0.6}{
        \tikzset{every picture/.style={line width=0.75pt}} 
        \begin{tikzpicture}[x=0.75pt,y=0.75pt,yscale=-1,xscale=1] 
        \draw  [color={rgb, 255:red, 0; green, 0; blue, 0 }  ,draw opacity=0 ][fill={rgb, 255:red, 155; green, 155; blue, 155 }  ,fill opacity=0.6 ] (300,80) -- (530,380) -- (300,380) -- cycle ;
        \draw  [color={rgb, 255:red, 0; green, 0; blue, 0 }  ,draw opacity=0 ][fill={rgb, 255:red, 155; green, 155; blue, 155 }  ,fill opacity=0.6 ] (300,80) -- (300,380) -- (80,380) -- cycle ;
        \draw [line width=2.25]    (300,80) -- (300,380) ;
        \draw [line width=2.25]    (530,380) -- (300,380) ;
        \draw  [fill={rgb, 255:red, 0; green, 0; blue, 0 }  ,fill opacity=1 ] (295,380) .. controls (295,377.24) and (297.24,375) .. (300,375) .. controls (302.76,375) and (305,377.24) .. (305,380) .. controls (305,382.76) and (302.76,385) .. (300,385) .. controls (297.24,385) and (295,382.76) .. (295,380) -- cycle ;
        \draw    (301.12,431.75) -- (319.23,451.15) ;
        \draw    (319.23,412.36) -- (301.12,431.75) ;
        \draw    (301.12,431.75) -- (301.12,405.89) ;
        \draw    (283,412.36) -- (301.12,431.75) ;
        \draw    (301.12,431.75) -- (283,451.15) ;
        \draw  [fill={rgb, 255:red, 0; green, 0; blue, 0 }  ,fill opacity=1 ] (298.1,431.75) .. controls (298.1,429.97) and (299.45,428.52) .. (301.12,428.52) .. controls (302.78,428.52) and (304.14,429.97) .. (304.14,431.75) .. controls (304.14,433.54) and (302.78,434.98) .. (301.12,434.98) .. controls (299.45,434.98) and (298.1,433.54) .. (298.1,431.75) -- cycle ;
        \draw    (106.23,221.06) -- (160.35,221.06) ;
        \draw  [fill={rgb, 255:red, 0; green, 0; blue, 0 }  ,fill opacity=1 ] (103.53,221.06) .. controls (103.53,219.31) and (104.74,217.9) .. (106.23,217.9) .. controls (107.73,217.9) and (108.94,219.31) .. (108.94,221.06) .. controls (108.94,222.8) and (107.73,224.21) .. (106.23,224.21) .. controls (104.74,224.21) and (103.53,222.8) .. (103.53,221.06) -- cycle ;
        \draw  [fill={rgb, 255:red, 0; green, 0; blue, 0 }  ,fill opacity=1 ] (130.59,221.06) .. controls (130.59,219.31) and (131.8,217.9) .. (133.29,217.9) .. controls (134.79,217.9) and (136,219.31) .. (136,221.06) .. controls (136,222.8) and (134.79,224.21) .. (133.29,224.21) .. controls (131.8,224.21) and (130.59,222.8) .. (130.59,221.06) -- cycle ;
        \draw  [fill={rgb, 255:red, 0; green, 0; blue, 0 }  ,fill opacity=1 ] (157.64,221.06) .. controls (157.64,219.31) and (158.86,217.9) .. (160.35,217.9) .. controls (161.84,217.9) and (163.06,219.31) .. (163.06,221.06) .. controls (163.06,222.8) and (161.84,224.21) .. (160.35,224.21) .. controls (158.86,224.21) and (157.64,222.8) .. (157.64,221.06) -- cycle ;
        \draw    (133.29,195.8) -- (133.29,221.06) ;
        \draw    (90,202.11) -- (106.23,221.06) ;
        \draw    (106.23,221.06) -- (90,240) ;
        \draw    (160.35,221.06) -- (176.59,240) ;
        \draw    (176.59,202.11) -- (160.35,221.06) ;
        \draw    (277.27,65.96) -- (309.38,65.96) ;
        \draw  [fill={rgb, 255:red, 0; green, 0; blue, 0 }  ,fill opacity=1 ] (274.06,65.96) .. controls (274.06,64.03) and (275.49,62.48) .. (277.27,62.48) .. controls (279.04,62.48) and (280.48,64.03) .. (280.48,65.96) .. controls (280.48,67.88) and (279.04,69.44) .. (277.27,69.44) .. controls (275.49,69.44) and (274.06,67.88) .. (274.06,65.96) -- cycle ;
        \draw  [fill={rgb, 255:red, 0; green, 0; blue, 0 }  ,fill opacity=1 ] (306.17,65.96) .. controls (306.17,64.03) and (307.61,62.48) .. (309.38,62.48) .. controls (311.15,62.48) and (312.59,64.03) .. (312.59,65.96) .. controls (312.59,67.88) and (311.15,69.44) .. (309.38,69.44) .. controls (307.61,69.44) and (306.17,67.88) .. (306.17,65.96) -- cycle ;
        \draw    (258,45.08) -- (277.27,65.96) ;
        \draw    (277.27,65.96) -- (258,86.83) ;
        \draw    (309.38,65.96) -- (328.65,86.83) ;
        \draw    (328.65,45.08) -- (309.38,65.96) ;
        \draw    (309.38,65.96) -- (328.65,65.96) ;
        \draw    (441.14,215.27) -- (471.82,215.27) ;
        \draw  [fill={rgb, 255:red, 0; green, 0; blue, 0 }  ,fill opacity=1 ] (438.07,215.27) .. controls (438.07,213.45) and (439.45,211.98) .. (441.14,211.98) .. controls (442.84,211.98) and (444.21,213.45) .. (444.21,215.27) .. controls (444.21,217.08) and (442.84,218.55) .. (441.14,218.55) .. controls (439.45,218.55) and (438.07,217.08) .. (438.07,215.27) -- cycle ;
        \draw  [fill={rgb, 255:red, 0; green, 0; blue, 0 }  ,fill opacity=1 ] (468.76,215.27) .. controls (468.76,213.45) and (470.13,211.98) .. (471.82,211.98) .. controls (473.52,211.98) and (474.89,213.45) .. (474.89,215.27) .. controls (474.89,217.08) and (473.52,218.55) .. (471.82,218.55) .. controls (470.13,218.55) and (468.76,217.08) .. (468.76,215.27) -- cycle ;
        \draw    (422.73,195.56) -- (441.14,215.27) ;
        \draw    (441.14,215.27) -- (422.73,234.98) ;
        \draw    (471.82,215.27) -- (490.23,234.98) ;
        \draw    (490.23,195.56) -- (471.82,215.27) ;
        \draw    (471.82,215.27) -- (490.23,215.27) ;
        \draw    (492.64,419.67) -- (527.61,419.67) ;
        \draw  [fill={rgb, 255:red, 0; green, 0; blue, 0 }  ,fill opacity=1 ] (489.14,419.67) .. controls (489.14,417.8) and (490.71,416.29) .. (492.64,416.29) .. controls (494.57,416.29) and (496.14,417.8) .. (496.14,419.67) .. controls (496.14,421.54) and (494.57,423.06) .. (492.64,423.06) .. controls (490.71,423.06) and (489.14,421.54) .. (489.14,419.67) -- cycle ;
        \draw  [fill={rgb, 255:red, 0; green, 0; blue, 0 }  ,fill opacity=1 ] (524.11,419.67) .. controls (524.11,417.8) and (525.67,416.29) .. (527.61,416.29) .. controls (529.54,416.29) and (531.1,417.8) .. (531.1,419.67) .. controls (531.1,421.54) and (529.54,423.06) .. (527.61,423.06) .. controls (525.67,423.06) and (524.11,421.54) .. (524.11,419.67) -- cycle ;
        \draw    (471.66,399.37) -- (492.64,419.67) ;
        \draw    (492.64,419.67) -- (471.66,439.98) ;
        \draw    (527.61,419.67) -- (548.59,439.98) ;
        \draw    (548.59,399.37) -- (527.61,419.67) ;
        \draw    (527.61,419.67) -- (548.59,419.67) ;
        
        \draw (324,29) node [anchor=north west][inner sep=0.75pt]    [xscale=1.5, yscale=1.5] {$\infty $};

        \draw (121,174) node [anchor=north west][inner sep=0.75pt]   [xscale=1.5, yscale=1.5] {$\infty $};
        
        \draw (406,175) node [anchor=north west][inner sep=0.75pt]    [xscale=1.5, yscale=1.5] {$\infty $};

        \draw (290,388.85) node [anchor=north west][inner sep=0.75pt] [xscale=1.5, yscale=1.5] {$\infty $};

        \draw (452, 383) node [anchor=north west][inner sep=0.75pt] [xscale=1.5, yscale=1.5] {$\infty $};
        
        \end{tikzpicture} }
    \caption{The space $M^{\mathrm{trop}}_{0,5}/S_{4}$}
    \label{fig:M_trop_041}
\end{figure}

In Figures~\ref{fig:M_trop_05} and \ref{fig:M_trop_041}, the (possibly folded) positive orthants are glued with respect to degeneration of the corresponding tree types, which is shown in Figure~\ref{fig:GraphDegenerations}. We refer to \cite[Section~2]{MarkwigSurvey} and \cite[Section~4]{ChanLectures} for more details. 

 \begin{figure}[ht]
    \centering

    \scalebox{0.6}{  
        \tikzset{every picture/.style={line width=0.75pt}}
        \begin{tikzpicture}[x=0.75pt,y=0.75pt,yscale=-1,xscale=1]
        \draw    (292.68,85.61) -- (378.05,85.61) ;
        \draw  [fill={rgb, 255:red, 0; green, 0; blue, 0 }  ,fill opacity=1 ] (288.41,85.61) .. controls (288.41,83.05) and (290.33,80.98) .. (292.68,80.98) .. controls (295.04,80.98) and (296.95,83.05) .. (296.95,85.61) .. controls (296.95,88.17) and (295.04,90.24) .. (292.68,90.24) .. controls (290.33,90.24) and (288.41,88.17) .. (288.41,85.61) -- cycle ;
        \draw  [fill={rgb, 255:red, 0; green, 0; blue, 0 }  ,fill opacity=1 ] (331.1,85.61) .. controls (331.1,83.05) and (333.01,80.98) .. (335.37,80.98) .. controls (337.72,80.98) and (339.63,83.05) .. (339.63,85.61) .. controls (339.63,88.17) and (337.72,90.24) .. (335.37,90.24) .. controls (333.01,90.24) and (331.1,88.17) .. (331.1,85.61) -- cycle ;
        \draw  [fill={rgb, 255:red, 0; green, 0; blue, 0 }  ,fill opacity=1 ] (373.78,85.61) .. controls (373.78,83.05) and (375.69,80.98) .. (378.05,80.98) .. controls (380.41,80.98) and (382.32,83.05) .. (382.32,85.61) .. controls (382.32,88.17) and (380.41,90.24) .. (378.05,90.24) .. controls (375.69,90.24) and (373.78,88.17) .. (373.78,85.61) -- cycle ;
        \draw    (335.37,48.54) -- (335.37,59.75) -- (335.37,85.61) ;
        \draw    (267.07,57.8) -- (292.68,85.61) ;
        \draw    (292.68,85.61) -- (267.07,113.41) ;
        \draw    (378.05,85.61) -- (403.66,113.41) ;
        \draw    (403.66,57.8) -- (378.05,85.61) ;
        \draw    (489.02,85.61) -- (574.39,85.61) ;
        \draw  [fill={rgb, 255:red, 0; green, 0; blue, 0 }  ,fill opacity=1 ] (484.76,85.61) .. controls (484.76,83.05) and (486.67,80.98) .. (489.02,80.98) .. controls (491.38,80.98) and (493.29,83.05) .. (493.29,85.61) .. controls (493.29,88.17) and (491.38,90.24) .. (489.02,90.24) .. controls (486.67,90.24) and (484.76,88.17) .. (484.76,85.61) -- cycle ;
        \draw  [fill={rgb, 255:red, 0; green, 0; blue, 0 }  ,fill opacity=1 ] (527.44,85.61) .. controls (527.44,83.05) and (529.35,80.98) .. (531.71,80.98) .. controls (534.06,80.98) and (535.98,83.05) .. (535.98,85.61) .. controls (535.98,88.17) and (534.06,90.24) .. (531.71,90.24) .. controls (529.35,90.24) and (527.44,88.17) .. (527.44,85.61) -- cycle ;
        \draw  [fill={rgb, 255:red, 0; green, 0; blue, 0 }  ,fill opacity=1 ] (570.12,85.61) .. controls (570.12,83.05) and (572.03,80.98) .. (574.39,80.98) .. controls (576.75,80.98) and (578.66,83.05) .. (578.66,85.61) .. controls (578.66,88.17) and (576.75,90.24) .. (574.39,90.24) .. controls (572.03,90.24) and (570.12,88.17) .. (570.12,85.61) -- cycle ;
        \draw    (531.71,48.54) -- (531.71,85.61) ;
        \draw    (463.41,57.8) -- (489.02,85.61) ;
        \draw    (489.02,85.61) -- (463.41,113.41) ;
        \draw    (574.39,85.61) -- (600,113.41) ;
        \draw    (600,57.8) -- (574.39,85.61) ;
        \draw    (335.37,131.95) -- (427.55,186.54) ;
        \draw [shift={(429.27,187.56)}, rotate = 210.63] [color={rgb, 255:red, 0; green, 0; blue, 0 }  ][line width=0.75]    (10.93,-3.29) .. controls (6.95,-1.4) and (3.31,-0.3) .. (0,0) .. controls (3.31,0.3) and (6.95,1.4) .. (10.93,3.29)   ;
        \draw    (335.37,131.95) -- (335.37,185.56) ;
        \draw [shift={(335.37,187.56)}, rotate = 270] [color={rgb, 255:red, 0; green, 0; blue, 0 }  ][line width=0.75]    (10.93,-3.29) .. controls (6.95,-1.4) and (3.31,-0.3) .. (0,0) .. controls (3.31,0.3) and (6.95,1.4) .. (10.93,3.29)   ;
        \draw    (506.1,243.17) -- (548.78,243.17) ;
        \draw  [fill={rgb, 255:red, 0; green, 0; blue, 0 }  ,fill opacity=1 ] (501.83,243.17) .. controls (501.83,240.61) and (503.74,238.54) .. (506.1,238.54) .. controls (508.45,238.54) and (510.37,240.61) .. (510.37,243.17) .. controls (510.37,245.73) and (508.45,247.8) .. (506.1,247.8) .. controls (503.74,247.8) and (501.83,245.73) .. (501.83,243.17) -- cycle ;
        \draw  [fill={rgb, 255:red, 0; green, 0; blue, 0 }  ,fill opacity=1 ] (544.51,243.17) .. controls (544.51,240.61) and (546.42,238.54) .. (548.78,238.54) .. controls (551.14,238.54) and (553.05,240.61) .. (553.05,243.17) .. controls (553.05,245.73) and (551.14,247.8) .. (548.78,247.8) .. controls (546.42,247.8) and (544.51,245.73) .. (544.51,243.17) -- cycle ;
        \draw    (480.49,215.37) -- (506.1,243.17) ;
        \draw    (506.1,243.17) -- (480.49,270.98) ;
        \draw    (548.78,243.17) -- (574.39,270.98) ;
        \draw    (574.39,215.37) -- (548.78,243.17) ;
        \draw    (548.78,243.17) -- (574.39,243.17) ;
        \draw    (531.71,131.95) -- (531.71,185.56) ;
        \draw [shift={(531.71,187.56)}, rotate = 270] [color={rgb, 255:red, 0; green, 0; blue, 0 }  ][line width=0.75]    (10.93,-3.29) .. controls (6.95,-1.4) and (3.31,-0.3) .. (0,0) .. controls (3.31,0.3) and (6.95,1.4) .. (10.93,3.29)   ;
        \draw    (318.29,243.17) -- (360.98,243.17) ;
        \draw  [fill={rgb, 255:red, 0; green, 0; blue, 0 }  ,fill opacity=1 ] (314.02,243.17) .. controls (314.02,240.61) and (315.94,238.54) .. (318.29,238.54) .. controls (320.65,238.54) and (322.56,240.61) .. (322.56,243.17) .. controls (322.56,245.73) and (320.65,247.8) .. (318.29,247.8) .. controls (315.94,247.8) and (314.02,245.73) .. (314.02,243.17) -- cycle ;
        \draw  [fill={rgb, 255:red, 0; green, 0; blue, 0 }  ,fill opacity=1 ] (356.71,243.17) .. controls (356.71,240.61) and (358.62,238.54) .. (360.98,238.54) .. controls (363.33,238.54) and (365.24,240.61) .. (365.24,243.17) .. controls (365.24,245.73) and (363.33,247.8) .. (360.98,247.8) .. controls (358.62,247.8) and (356.71,245.73) .. (356.71,243.17) -- cycle ;
        \draw    (292.68,215.37) -- (318.29,243.17) ;
        \draw    (318.29,243.17) -- (292.68,270.98) ;
        \draw    (360.98,243.17) -- (386.59,270.98) ;
        \draw    (386.59,215.37) -- (360.98,243.17) ;
        \draw    (360.98,243.17) -- (386.59,243.17) ;
        \draw    (335.37,270.98) -- (402.11,325.32) ;
        \draw [shift={(403.66,326.59)}, rotate = 219.16] [color={rgb, 255:red, 0; green, 0; blue, 0 }  ][line width=0.75]    (10.93,-3.29) .. controls (6.95,-1.4) and (3.31,-0.3) .. (0,0) .. controls (3.31,0.3) and (6.95,1.4) .. (10.93,3.29)   ;
        \draw    (523.17,270.98) -- (473.31,325.11) ;
        \draw [shift={(471.95,326.59)}, rotate = 312.65] [color={rgb, 255:red, 0; green, 0; blue, 0 }  ][line width=0.75]    (10.93,-3.29) .. controls (6.95,-1.4) and (3.31,-0.3) .. (0,0) .. controls (3.31,0.3) and (6.95,1.4) .. (10.93,3.29)   ;
        \draw    (437.8,382.2) -- (463.41,410) ;
        \draw    (463.41,354.39) -- (437.8,382.2) ;
        \draw    (437.8,382.2) -- (437.8,345.12) ;
        \draw    (412.2,354.39) -- (437.8,382.2) ;
        \draw    (437.8,382.2) -- (412.2,410) ;
        \draw  [fill={rgb, 255:red, 0; green, 0; blue, 0 }  ,fill opacity=1 ] (433.54,382.2) .. controls (433.54,379.64) and (435.45,377.56) .. (437.8,377.56) .. controls (440.16,377.56) and (442.07,379.64) .. (442.07,382.2) .. controls (442.07,384.75) and (440.16,386.83) .. (437.8,386.83) .. controls (435.45,386.83) and (433.54,384.75) .. (433.54,382.2) -- cycle ;
        \draw    (75.61,97.07) -- (160.98,97.07) ;
        \draw  [fill={rgb, 255:red, 0; green, 0; blue, 0 }  ,fill opacity=1 ] (71.34,97.07) .. controls (71.34,94.51) and (73.25,92.44) .. (75.61,92.44) .. controls (77.97,92.44) and (79.88,94.51) .. (79.88,97.07) .. controls (79.88,99.63) and (77.97,101.71) .. (75.61,101.71) .. controls (73.25,101.71) and (71.34,99.63) .. (71.34,97.07) -- cycle ;
        \draw  [fill={rgb, 255:red, 0; green, 0; blue, 0 }  ,fill opacity=1 ] (114.02,97.07) .. controls (114.02,94.51) and (115.94,92.44) .. (118.29,92.44) .. controls (120.65,92.44) and (122.56,94.51) .. (122.56,97.07) .. controls (122.56,99.63) and (120.65,101.71) .. (118.29,101.71) .. controls (115.94,101.71) and (114.02,99.63) .. (114.02,97.07) -- cycle ;
        \draw  [fill={rgb, 255:red, 0; green, 0; blue, 0 }  ,fill opacity=1 ] (156.71,97.07) .. controls (156.71,94.51) and (158.62,92.44) .. (160.98,92.44) .. controls (163.33,92.44) and (165.24,94.51) .. (165.24,97.07) .. controls (165.24,99.63) and (163.33,101.71) .. (160.98,101.71) .. controls (158.62,101.71) and (156.71,99.63) .. (156.71,97.07) -- cycle ;
        \draw    (118.29,60) -- (118.29,71.21) -- (118.29,97.07) ;
        \draw    (50,69.27) -- (75.61,97.07) ;
        \draw    (75.61,97.07) -- (50,124.88) ;
        \draw    (160.98,97.07) -- (186.59,124.88) ;
        \draw    (186.59,69.27) -- (160.98,97.07) ;
        \draw    (120,134.39) -- (120,188) ;
        \draw [shift={(120,190)}, rotate = 270] [color={rgb, 255:red, 0; green, 0; blue, 0 }  ][line width=0.75]    (10.93,-3.29) .. controls (6.95,-1.4) and (3.31,-0.3) .. (0,0) .. controls (3.31,0.3) and (6.95,1.4) .. (10.93,3.29)   ;
        \draw    (95.61,232.2) -- (138.29,232.2) ;
        \draw  [fill={rgb, 255:red, 0; green, 0; blue, 0 }  ,fill opacity=1 ] (91.34,232.2) .. controls (91.34,229.64) and (93.25,227.56) .. (95.61,227.56) .. controls (97.97,227.56) and (99.88,229.64) .. (99.88,232.2) .. controls (99.88,234.75) and (97.97,236.83) .. (95.61,236.83) .. controls (93.25,236.83) and (91.34,234.75) .. (91.34,232.2) -- cycle ;
        \draw  [fill={rgb, 255:red, 0; green, 0; blue, 0 }  ,fill opacity=1 ] (134.02,232.2) .. controls (134.02,229.64) and (135.94,227.56) .. (138.29,227.56) .. controls (140.65,227.56) and (142.56,229.64) .. (142.56,232.2) .. controls (142.56,234.75) and (140.65,236.83) .. (138.29,236.83) .. controls (135.94,236.83) and (134.02,234.75) .. (134.02,232.2) -- cycle ;
        \draw    (70,204.39) -- (95.61,232.2) ;
        \draw    (95.61,232.2) -- (70,260) ;
        \draw    (138.29,232.2) -- (163.9,260) ;
        \draw    (163.9,204.39) -- (138.29,232.2) ;
        \draw    (138.29,232.2) -- (163.9,232.2) ;
        \draw    (120,260) -- (120,328) ;
        \draw [shift={(120,330)}, rotate = 270] [color={rgb, 255:red, 0; green, 0; blue, 0 }  ][line width=0.75]    (10.93,-3.29) .. controls (6.95,-1.4) and (3.31,-0.3) .. (0,0) .. controls (3.31,0.3) and (6.95,1.4) .. (10.93,3.29)   ;
        \draw    (119.61,382.2) -- (145.22,410) ;
        \draw    (145.22,354.39) -- (119.61,382.2) ;
        \draw    (119.61,382.2) -- (119.61,345.12) ;
        \draw    (94,354.39) -- (119.61,382.2) ;
        \draw    (119.61,382.2) -- (94,410) ;
        \draw  [fill={rgb, 255:red, 0; green, 0; blue, 0 }  ,fill opacity=1 ] (115.34,382.2) .. controls (115.34,379.64) and (117.25,377.56) .. (119.61,377.56) .. controls (121.97,377.56) and (123.88,379.64) .. (123.88,382.2) .. controls (123.88,384.75) and (121.97,386.83) .. (119.61,386.83) .. controls (117.25,386.83) and (115.34,384.75) .. (115.34,382.2) -- cycle ;
        \draw (245,30) node [anchor=north west][inner sep=0.75pt]  [xscale=1.5, yscale=1.5] {$\infty $};
        \draw (516, 30) node [anchor=north west][inner sep=0.75pt]     [xscale=1.5, yscale=1.5] {$\infty $};
        \draw (574, 200) node [anchor=north west][inner sep=0.75pt]     [xscale=1.5, yscale=1.5] {$\infty $};
        \draw (270, 200) node [anchor=north west][inner sep=0.75pt]      [xscale=1.5, yscale=1.5] {$\infty $};
        \draw (426, 324) node [anchor=north west][inner sep=0.75pt]     [xscale=1.5, yscale=1.5] {$\infty $};
    \end{tikzpicture}

    }
    \caption{Degenerations of trees}
    \label{fig:GraphDegenerations}
\end{figure}
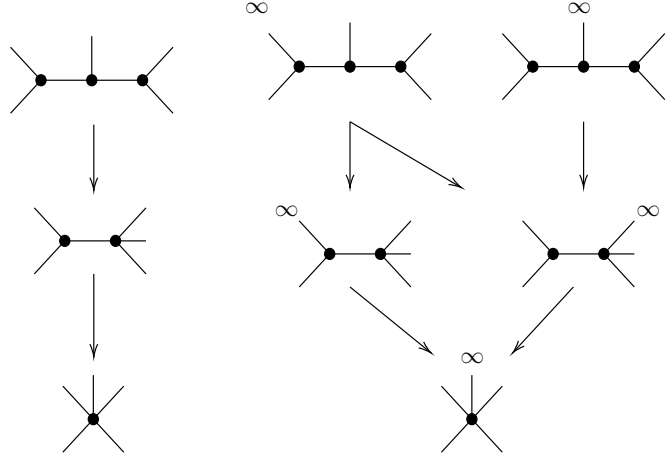

 \end{example}

For the symmetrized moduli space $M_{0,n}/S_{\Xi}$, we can again introduce an abstract tropicalization map as follows. Let $\Xi\subset \mathcal{P}[n]$ be a partition and let $P$ be an $L$-valued point of $M^{\mathrm{an}}_{0,n}/S_{\Xi}$. This again gives a phylogenetic tree $(T,\ell)$ with a marking $\{1,...,n\}\to L(T)$, but we identify two markings $\phi_{1},\phi_{2}:\{1,...,n\}\to L(T)$ if they are related by an element of $S_{\Xi}$. This gives a map
 \begin{equation}\label{TropSn}     \mathrm{trop}:M^{\mathrm{an}}_{0,n}/S_{\Xi}\to{M^{\mathrm{trop}}_{0,n}/S_{\Xi}}.
 \end{equation}
 
 \begin{remark}\label{CanDecomOfTropModSp} 
  The map in \eqref{TropSn} in particular gives a partition of $M^{\mathrm{an}}_{0,n}/S_{\Xi}$ by considering the inverse images of the loci of $\Xi$-marked tree types.
  We will use this partition to obtain a partition of the space of binary forms in Section~\ref{sec:TreesTropicalBinaryForms}. We note that the partition of $M^{\mathrm{an}}_{0,n}/S_{\Xi}$ is a special instance of the partitions we considered in Section~\ref{sec:CompactificationsTropicalizations} since the construction in \eqref{eq:ValuativeLimit} gives the reduction map from \eqref{eq:ReductionMap}. 
 \end{remark}
 
    \subsection{Trees and tropical binary forms}\label{sec:TreesTropicalBinaryForms}

    Let $V=V_{k_{1}}\oplus \cdots \oplus{V_{k_{r}}}$ be the standard $\SL_{2}$-module and write $A^{G}$ for the ring of invariants. We consider $\mathrm{Proj}(A^{G})$, where the grading is induced by the degree of an invariant. We are interested in the affine open $U=D_{+}(\Delta)$ for $\Delta$ the discriminant. The set of $L$-valued points $U(L)$ of $U$ can be identified with tuples of binary forms $(f_{1},\dots,f_{r})$ over $L$ without common zeros and up to $\GL_{2}$-equivalence. Here $\mathrm{deg}(f_{i})=k_{i}$. 
    
    We now consider an $L$-valued point of $U$. This gives a set of binary forms $(f_{1},\dots,f_{r})$ over $L$. We write $f=\prod_{i=1}^{r}f_{i}$ for the resulting $(k_{1} + \cdots + k_{r})$-binary form. Let $Z(f)$ and $Z(f_i)$ denote the zero sets of $f$ and $f_i$ respectively, for $1 \leq i \leq r$. Note that $Z(f) = \sqcup_{i=1}^{r} \ Z(f_i)$ gives a well-defined $L$-valued point of the moduli space ${M}_{0,n}/S_{\Xi}$, where $n = k_1 + \cdots + k_r$ and $\Xi = \{ \{1, \dots, k_1\}, \{k_1 +1, \dots, k_1 + k_2\} , \dots , \{ k_1 + \cdots + k_{r-1} + 1, \dots , n\} \}$. Indeed, taking a different equivalent binary form changes the zero sets by the action of $\GL_{2}$, so that the induced map 
    \begin{equation*}
        \left(\mathbb{P}^{1}, \ \bigsqcup_{i=1}^{r} Z(f_i)\right) \to \left(\mathbb{P}^{1}, \ \bigsqcup_{i=1}^{r} Z(f_i^\sigma)\right)
    \end{equation*} is an isomorphism, which means that we obtain the same point in ${M}_{0,n}/S_{\Xi}$. 
    \begin{definition}
        \textbf{[Metric tree of a binary form]}
        Let $f=(f_{1},\dots,f_{r})$ be a separable binary form over a valued field $L$ and consider a set of markings of the $Z(f_{i})$. Let $\{1,...,n\}\to Z(f)$ be the induced marking of $Z(f)$ with marked metric tree $(T,\ell)$, see \cref{def:MarkedMetricTree}. Then the \emph{metric tree} associated to $f$ is the image of $(T,\ell)$ in $M^{\mathrm{trop}}_{0,n}/S_{\Xi}$ under the map $M^{\mathrm{trop}}_{0,n}\to M^{\mathrm{trop}}_{0,n}/S_{\Xi}$. For any $\sigma\in\GL_{2}$ with binary form $\sigma(f)=(\sigma(f_{1}),...,\sigma(f_{r}))$, the induced point of $M^{\mathrm{trop}}_{0,n}/S_{\Xi}$ is the same. 
    \end{definition}
    The partition on the Berkovich analytification of $M_{0,n}/S_{\Xi}$ given in Remark~\ref{CanDecomOfTropModSp} induces a partition of $U^{\mathrm{an}}$. We use this as our definition of a tropical invariant of a binary form. 

    \begin{definition}
    {\bf{[Tropical invariants of binary forms]}}
    A set of \textit{tropical invariants} of a set of binary forms of degree $n$ is a tropical separating set for the partition of $U^{\an}$ induced from the maps $U^{\an} \to M^{\mathrm{an}}_{0,n}/S_{\Xi}\to{M^{\mathrm{trop}}_{0,n}/S_{\Xi}}$, see \eqref{TropSn}.
    \end{definition}

    For a given set of invariants $h_{i}$, we obtain a rational map $U\to\mathbb{P}(a_{0},\dots,a_{n})$,
    where $\mathrm{deg}(h_{i})=a_{i}$. If the $h_{i}$ include a set of generators of the ring of invariants, then this is automatically a morphism since the set of generators generate the nullcone. If the set of generators is a projective tropical separating set, then we call this a set of \textit{projective tropical invariants}.

    \begin{example}
    For residue characteristics not equal to $2$, the tropical invariants of a binary quartic are given by $c_{4}$, $c_{6}$ and $\Delta$; see \cite[Chapter~VII]{SilvermanI}. These are furthermore weighted projective invariants, giving a tropicalization map to $\mathbb{T}\mathbb{P}(4,6,12)$. 
    \end{example}
    
    \begin{example}
    For residue characteristics not equal to $2$, the tropical invariants of a binary sextic are given by the Igusa invariants from \cite{helminck2021Igusa}.  
    \end{example}
    
    \begin{example}
    In Section~\ref{Sec:ProofMainResults}, we will show that the invariants introduced in the introduction form a set of tropical invariants for binary quintics and $(4,1)$-forms. These then also give the reduction types of Picard curves when the residue characteristic is not $3$.  
    \end{example}

    As promised in Section~\ref{subsection:quinticInvariants}, we now show that the invariants $I_{4}$, $I_{8}$, $I_{12}$, $I_{18}$ and $\Delta$ are not enough to distinguish between the unmarked trees of a quintic.

	\begin{example}\label{ex:counterExampleTreeType}
	   Let $K=\mathbb{C}\{\{t\}\}$ be the field of Puiseux series over the complex numbers with $v(t)=1$. Consider the two quintics given by
	   \begin{align*}
	       f_{2}& = xz(x-z)(x-t^2z)(x-2z),\\
	       f_{3}& = xz(x-z)(x-tz)(x-(1+2t)z).
	   \end{align*}
	   The tree of $f_{2}$ is of Type~$\mathrm{II}$ and the tree of $f_{3}$ is of Type~$\mathrm{III}$. However, in both cases, the tropical invariants are the same:
	   \[
	    [v(I_{4}) : v(I_{8}) : v(I_{12}) : v(I_{18}) : v(\Delta)] = [0:0:0:2:4].
	   \]
	   We thus see that we cannot distinguish between the two tree types using these invariants. This also implies that we cannot distinguish the various reduction types of Picard curves using these invariants. They are however enough to decide whether a tree is of Type~$\mathrm{I}$ or not.

	\end{example}

    The results in this paper show that there exists a set of projective tropical invariants for quintics and $(4,1)$-forms. In general, we conjecture that there exists a finite set of projective tropical invariants for any set of binary forms. Moreover, there should be a practical algorithm that can calculate these tropical invariants.

    \section{Proofs of the main results} \label{Sec:ProofMainResults}
    
        In this final section, we prove the main results stated in the introduction, namely \cref{thm:tropicalTypes}, \cref{thm:reductionTypes}, \cref{thm:edgeLengths} and Corollary~\ref{cor:reductionTypesOfPicardCurves}. Some parts of the proofs rely on computing the invariants explicitly. The computations are made available in (\ref{eq:link_for_code}).

        Recall that our base field $K$ is a non-archimedean, complete and algebraically closed valued field of characteristic $0$ with associated data $(v,R,\mathfrak{m},k)$. Moreover, the residue characteristic $p$ is different from $2,3$ and $11$.
        
        \begin{remark} \label{rem:Thm1p=11}
                If the residue characteristic of $K$ is $11$, then the second condition in \cref{thm:tropicalTypes} does not characterize trees of Type~$\mathrm{II}$. In this case, to obtain condition for Type~$\mathrm{II}$, we simply need to check that the conditions for Type~$\mathrm{I}$ and Type~$\mathrm{III}$ are not satisfied.
        \end{remark}

        \subsection{Universal families}\label{sec:UniversalFamilies}

        In this section, we explain the general idea of the proofs. We begin by recalling the setup in Sections \ref{sec:ModuliStable} and \ref{sec:TreesTropicalBinaryForms} in the case of binary quintics and $(4,1)$-forms.
        
        Let $f$ be a separable binary quintic or $(4,1)$-form over $K$. The zero set $Z(f)$ of $f$ consists of five points in $\mathbb{P}^{1}$, which we also consider as points of the corresponding Berkovich analytification $\mathbb{P}^{1,\mathrm{an}}$.
        These are often called points of type $1$ in the literature. These five points define a natural metric tree $\Sigma_{f}\subset\mathbb{P}^{1,\an}$ that contains $Z(f)$. If we choose a marking $\{1,...,5\}\to Z(f)$, then this is the \emph{minimal Berkovich skeleton} of the marked curve $(\mathbb{P}^{1},P_{1},...,P_{5})$ as defined in \cite[Definition 4.20]{Baker_Payne_Rabinoff_13} in terms of semistable vertex sets. If $f=(f_{4},f_{1})$ is a $(4,1)$-form, then we assume that the marking is chosen so that $P_{1},...,P_{4}\in Z(f_{4})$.

        The group $\GL_{2}(K)$ acts on $\mathbb{P}^{1,\an}$ through M\"{o}bius transformations. 
        The individual $\Sigma_{f}$ are not invariant under this action, since $Z(f)$ is not invariant. The abstract marked metric tree $\Sigma_{f}$ however is invariant, up to suitable permutations of the markings. More precisely, if $f$ is a binary quintic, then we consider $\Sigma_{f}$ as an element of the moduli space $M^{\mathrm{trop}}_{0,5}/S_{5}$ of unmarked metric trees on five leaves, and if $f=(f_{4},f_{1})$ is a $(4,1)$-form, then we consider $\Sigma_{f}$ as an element of the moduli space $M^{\mathrm{trop}}_{0,5}/S_{4}$ of metric trees with five leaves and one marked point. For every separable binary form in the same $\mathrm{GL}_{2}$-orbit, this gives the same metric tree, so that we can associate an element of $M^{\mathrm{trop}}_{0,5}/S_{5}$ or $M^{\mathrm{trop}}_{0,5}/S_{4}$ to any $\mathrm{GL}_{2}$-equivalence class of separable binary quintics or $(4,1)$-forms. Note that if two separable binary forms lie in the same $\mathrm{GL}_{2}$-orbit, then they also have the same projective invariants since a transformation $\sigma$ scales an invariant $I$ by
            $\mathrm{det}(\sigma)^{\mathrm{deg}(I)}$; see Remark~\ref{rem:GL2vsSL2}. To summarize, applying a transformation in $\mathrm{GL}_{2}$ to a binary quintic or $(4,1)$-form does not change
            \begin{enumerate}
                \item the metric tree of the binary form,
                \item the projective invariants of the binary form.
            \end{enumerate}

            \medskip
            
               For a given binary quintic or $(4,1)$-form defined over a valued field, we can now apply a projective transformation so that the resulting binary form is 
        \begin{equation}\label{eq:BinaryQuinticStandard}
            f(x,z) = xz(x-z)(x-\lambda_{1}z)(x-\lambda_{2}z).
        \end{equation}
         Here, we send the marked point to $\infty$.
        By the discussion above, this does not change the associated metric tree, nor the projective invariants.
        We now view the binary form in \eqref{eq:BinaryQuinticStandard} as being defined over the ring $A_{0}=K[\lambda_{1},\lambda_{2}]$. The latter is the coordinate ring of a two-dimensional affine space $Y_{0}=\mathbb{A}^{2}_{K}$ over $K$. For specializations of the $\lambda_{i}$ outside the vanishing locus $Z(\Delta_{\lambda})$ of $ \Delta_{\lambda} \coloneqq \lambda_{1}\lambda_{2}(\lambda_{1}-1)(\lambda_{2}-1)(\lambda_{1}-\lambda_{2})$, we again obtain a separable binary form. We write $Y$ for this open subspace of $Y_{0}$ and $U$ for the space of separable binary forms as in \cref{sec:TreesTropicalBinaryForms}. We then have a natural map 
        $$Y\to U$$ 
        sending $(\lambda_{1},\lambda_{2})$ to the separable binary quintic in \eqref{eq:BinaryQuinticStandard}. This map is surjective by construction. 

        We will write down parametrizations of the strata of $Y^{\an}$ that correspond to the different tree types in $U^{\an}$. This will be done in terms of \emph{universal algebras}. For each tree type, we will give a free algebra $A=R[t_{i},\mu_{i}]$ with an injective map $R[\lambda_{1},\lambda_{2}] \to  A$. We will view this as an inclusion, so that $\lambda_{i}\in{A}$. Note that the binary form in \eqref{eq:BinaryQuinticStandard} automatically becomes a binary form over $A$ through this map. Let $R_{L}$ be the valuation ring of a valued field extension $K\subset{L}$. We then consider  $R_{L}$-valued points $\psi:A\to{R_{L}}$ with $v(\psi(t_{i}))>0$. Such an $R_{L}$-valued point corresponds to a specific binary form of the given type. We will also write $v(t_{i})>0$ if $\psi$ is understood. The universal algebras are constructed in such a way that every separable binary form of the given type occurs as an $R_{L}$-valued point with $v(t_{i})>0$.

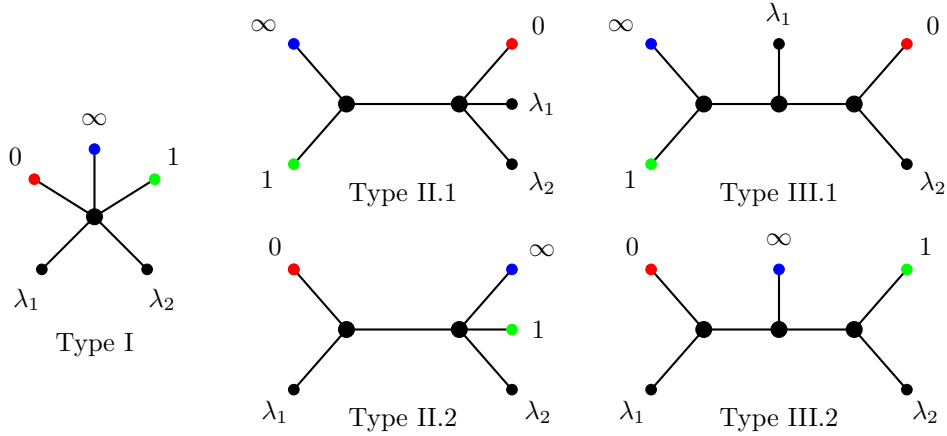
\begin{figure}[H]
    	    \centering
            	    \begin{tikzpicture}
                        \filldraw (0.9,0) circle (3pt);
                        \draw [thick] (0.9,0) -- (0.9,0.9);
                        \draw [thick] (0.9,0) -- (1.7,0.5);
                        \draw [thick] (0.9,0) -- (0.1,0.5);
                        \draw [thick] (0.9,0) -- (1.6,-0.7);
                        \draw [thick] (0.9,0) -- (0.2,-0.7);
                        
                        \filldraw[red] (0.1,0.5) circle (2pt);
                        \filldraw (-0.15,0.8) node{$0$};
                        
                        \filldraw[blue] (0.9,0.9) circle (2pt);
                        \filldraw (0.9,1.3) node{$\infty$};
                        
                        \filldraw[green] (1.7,0.5) circle (2pt);
                        \filldraw (1.95,0.8) node{$1$};
                        
                        \filldraw (0.2,-0.7) circle (2pt);
                        \filldraw (0,-1.1) node{$\lambda_1$};
                        
                        \filldraw (1.6,-0.7) circle (2pt);
                        \filldraw (1.8,-1.1) node{$\lambda_2$};
                        
                        \filldraw (0.9,-1.7) node{Type~$\mathrm{I}$};
                        
                        \qquad
                        
                        \filldraw (4.25,1.5) circle (3pt);
                        \filldraw (5.75,1.5) circle (3pt);
                        \draw [thick] (4.25,1.5) -- (5.75,1.5);
                        \draw [thick] (4.25,1.5) -- (3.55,2.3);
                        \draw [thick] (4.25,1.5) -- (3.55,0.7);
                        \draw [thick] (5.75,1.5) -- (6.45,2.3);
                        \draw [thick] (5.75,1.5) -- (6.45,1.5);
                        \draw [thick] (5.75,1.5) -- (6.45,0.7);
                        
                        \filldraw[blue] (3.55,2.3) circle (2pt);
                        \filldraw (3.15,2.55) node{$\infty$};
                        
                        \filldraw[green] (3.55,0.7) circle (2pt);
                        \filldraw (3.2,0.5) node{$1$};
                        
                        \filldraw[red] (6.45,2.3) circle (2pt);
                        \filldraw (6.8,2.55) node{$0$};
                        
                        \filldraw (6.45,1.5) circle (2pt);
                        \filldraw (6.85,1.5) node{$\lambda_1$};
                        
                        \filldraw (6.45,0.7) circle (2pt);
                        \filldraw (6.85,0.45) node{$\lambda_2$};
                
                        \filldraw (5,0.3) node{Type~$\mathrm{II}.1$};
                        
                        \filldraw (4.25,-1.5) circle (3pt);
                        \filldraw (5.75,-1.5) circle (3pt);
                        \draw [thick] (4.25,-1.5) -- (5.75,-1.5);
                        \draw [thick] (4.25,-1.5) -- (3.55,-0.7);
                        \draw [thick] (4.25,-1.5) -- (3.55,-2.3);
                        \draw [thick] (5.75,-1.5) -- (6.45,-0.7);
                        \draw [thick] (5.75,-1.5) -- (6.45,-1.5);
                        \draw [thick] (5.75,-1.5) -- (6.45,-2.3);
                        
                        \filldraw[red] (3.55,-0.7) circle (2pt);
                        
                        \filldraw[blue] (6.45,-0.7) circle (2pt);
                        \filldraw (6.85,-0.45) node{$\infty$};
                        
                        \filldraw[red] (3.55,-0.7) circle (2pt);
                        \filldraw (3.3,-0.4) node{$0$};
                        
                        \filldraw(3.55,-2.3) circle (2pt);
                        \filldraw (3.3,-2.6) node{$\lambda_1$};
                        
                        \filldraw[green] (6.45,-1.5) circle (2pt);
                        \filldraw (6.8,-1.5) node{$1$};
                        
                        \filldraw (6.45,-2.3) circle (2pt);
                        \filldraw (6.8,-2.6) node{$\lambda_2$};

                        \filldraw (5,-2.7) node{Type~$\mathrm{II}.2$};
                
                         \qquad
                         
                        \filldraw (9,1.5) circle (3pt);
                        \filldraw (10,1.5) circle (3pt);
                        \filldraw (11,1.5) circle (3pt);
                        \draw [thick] (9,1.5) -- (11,1.5);
                        \draw [thick] (9,1.5) -- (8.3,2.3);
                        \draw [thick] (9,1.5) -- (8.3,0.7);
                        \draw [thick] (10,1.5) -- (10,2.3);
                        \draw [thick] (11,1.5) -- (11.7,2.3);
                        \draw [thick] (11,1.5) -- (11.7,0.7);
                        
                        \filldraw[blue] (8.3,2.3) circle (2pt);
                        \filldraw (7.9,2.55) node{$\infty$};
                        
                        \filldraw[green] (8.3,0.7) circle (2pt);
                        \filldraw (8.025,0.5) node{$1$};
                        
                        \filldraw (10,2.3) circle (2pt);
                        \filldraw (10,2.7) node{$\lambda_1$};
                        
                        \filldraw[red] (11.7,2.3) circle (2pt);
                        \filldraw (12.05,2.55) node{$0$};
                        
                        \filldraw (11.7,0.7) circle (2pt);
                        \filldraw (12.05,0.45) node{$\lambda_2$};
                
                        \filldraw (10,0.3) node{Type~$\mathrm{III}.1$};

                        \filldraw (9,-1.5) circle (3pt);
                        \filldraw (10,-1.5) circle (3pt);
                        \filldraw (11,-1.5) circle (3pt);
                        \draw [thick] (9,-1.5) -- (11,-1.5);
                        \draw [thick] (9,-1.5) -- (8.3,-0.7);
                        \draw [thick] (9,-1.5) -- (8.3,-2.3);
                        \draw [thick] (10,-1.5) -- (10,-0.7);
                        \draw [thick] (11,-1.5) -- (11.7,-0.7);
                        \draw [thick] (11,-1.5) -- (11.7,-2.3);
                        
                        \filldraw[red] (8.3,-0.7) circle (2pt);
                        \filldraw (8.05,-0.4) node{$0$};
                        
                        \filldraw (8.3,-2.3) circle (2pt);
                        \filldraw (8.05,-2.6) node{$\lambda_1$};

                        \filldraw[blue] (10,-0.7) circle (2pt);
                        \filldraw (10,-0.3) node{$\infty$};
                        
                        \filldraw[green] (11.7,-0.7) circle (2pt);
                        \filldraw (11.95,-0.4) node{$1$};
                        
                        \filldraw (11.7,-2.3) circle (2pt);
                        \filldraw (11.95,-2.6) node{$\lambda_2$};
                
                        \filldraw (10,-2.7) node{Type~$\mathrm{III}.2$};
                
                        \end{tikzpicture}
    	    \caption{Trees corresponding to the universal families in Table~\ref{tab:UniFam}}
    	    \label{fig:TreesOfUnivFam}
    	\end{figure}

        To explicitly construct these universal families, we fix a tree type and write down the implied conditions on the valuations of the $\lambda_{i}$. We then parametrize these conditions using additional parameters $t_{i}$ and $\mu_{i}$. Our choices for these universal families can be found in Table~\ref{tab:UniFam}; see Figure~\ref{fig:TreesOfUnivFam} for the corresponding trees.
        
        \begin{example}
        To give an example, let us explain the universal family for Type III.2 in more detail. By assumption, we have the conditions
        \[
            v(\lambda_{1}) > 0 \quad \text{ and } \quad   v(\lambda_{2}-1) > 0.
        \]
        Note that we can assume by permuting the $\lambda_{i}$ that $v(\lambda_{1})\leq v(\lambda_{2}-1)$. To parametrize the different edge lengths, we then write 
        \[
            \lambda_{1}=t_{1}\mu_{1}  \quad \text{ and } \quad  \lambda_{2}=1+t_{2}\mu_{2}.
        \]
        Here we specialize the $\mu_{i}$ to elements with valuation zero and the $t_{i}$ to elements in $R$, so that $v(t_{1})=v(\lambda_{1})$ and $v(t_{2})=v(\lambda_{2}-1)$. In particular, note that the two non-trivial edge lengths are given exactly by $v(t_{i})$. By our assumption on the edge lengths, we then have $v(t_{1})\leq v(t_{2})$.   
        \end{example}

        \begin{table}[H]
            \centering
            \begin{tabular}{|c|c|c|c|}
                \hline
                \textbf{Type} & $\bm{\lambda_1}$ & $\bm{\lambda_2}$ & \textbf{Conditions}\\[0.6ex]
                \hline
                $\mathrm{I}$ & $\mu_1$ & $\mu_2$ & $\overline{\mu}_{i}\neq{0,1}$ and $\overline{\mu}_{1}\neq\overline{\mu}_{2}$ \\[0.6ex]
                \hline
                $\mathrm{II}.1$ &  $t_1\mu_1$ & $t_1\mu_2$  & $v(t_{1})>0$,  $\overline{\mu}_{i}\neq{0}$ and $\overline{\mu}_{1}\neq\overline{\mu}_{2}$ \\[0.6ex]
                \hline
                $\mathrm{II}.2$ &  $t_1\mu_1$ & $\mu_2$ & $v(t_{1})>0$, $\overline{\mu}_{i}\neq{0}$ and $\overline{\mu}_{2}\neq{1}$ \\[0.6ex]
                \hline
                $\mathrm{III}.1$ &  $t_1\mu_1$ & $t_1t_2\mu_2$ & $v(t_{i})>0$ and $\overline{\mu}_{i}\neq{0}$   \\[0.6ex]
                \hline
                $\mathrm{III}.2$ &  $t_1\mu_1$ & $1 + t_2\mu_2$ & $v(t_{i})>0$, $\overline{\mu}_{i}\neq{0}$ and $v(t_{1})\leq{v(t_{2})}$  \\[0.6ex]
                \hline
            \end{tabular}
        \caption{The chosen universal families. The $t_{i}$ parametrize the different edge lengths and the $\mu_{i}$ the different algebraic residue classes.}
        \label{tab:UniFam}
    \end{table}

        In these universal algebras, we have various natural reduction maps. We will be particularly interested in reducing elements modulo $I=\mathfrak{m}+(t_{1},t_{2})$. Here we think of the $t_{i}$
        as elements of $\mathfrak{m}$ whose valuations we can freely alter. In this vein, we will also refer to working modulo $I$ as working modulo $\mathfrak{m}$.
        For a given tree type with universal algebra $A$, we can consider the universal binary form $f$ from \eqref{eq:BinaryQuinticStandard} as being defined over $A$. The invariants of this binary form are thus elements of $A$.

        In particular, we now see that for any specialization $\psi:A\to{R_{L}}$, we have that $\psi(I)\in{R_{L}}$. We write $I\in{A^{\times}}$ if for every specialization $\psi$ as above, we have that $v(\psi(I))=0$.

      \subsection{Proof of \texorpdfstring{\cref{thm:tropicalTypes}}~}
      \label{pf:Thm1}
        
        In this subsection, we prove Theorem~\ref{thm:tropicalTypes}. In the proofs, we are free to choose the universal family of either $\mathrm{II}.1$ or $\mathrm{II}.2$, and similarly for $\mathrm{III}.1$ and $\mathrm{III}.2$. We do the necessity of the statement first and finish with the sufficiency in Section~\ref{sec:SufficiencyTheorem1}.
        
        \subsubsection{Type I}
        
        Computing the reduction modulo $\mathfrak{m}$ of the discriminant $\Delta$, we obtain
		    \[
		    \overline{\Delta} =  \overline{\lambda}_1^2 \  \overline{\lambda}_2^2 \ (\overline{\lambda}_1 - 1)^2 \ (\overline{\lambda}_2 - 1)^2 \ ( \overline{\lambda}_1 -\overline{\lambda}_2 )^2.
		    \]
		    So we deduce that $v(\Delta) = 0$, hence the condition $8 v(I) - \deg(I)v(\Delta) \geq 0$ is satisfied for any $I \in S$.

		    \subsubsection{Type II}

		    Consider the universal family for Type~$\mathrm{II}.1$ in Table~\ref{tab:UniFam}. We calculate the invariants $I \in S$ and find that they are divisible by $t_{1}^{\mathrm{deg}{I}/2}$. The reduction of $H/t_{1}^{\mathrm{deg}(H)/2}$ modulo $\mathfrak{m}$ is 
		    \begin{equation*}
		        \overline{H/t_{1}^{\mathrm{deg}(H)/2}}=-22\overline{\mu}_{1}^2\overline{\mu}_{2}^2(\overline{\mu}_{1}-\overline{\mu}_{2})^2.
		    \end{equation*}
		    Since $p\neq{2,11}$, we find that the latter is non-zero. We thus obtain 
		    \[
		    v(H) = \deg(H) v(t_1) /2 = 6 v(t_1), \quad \text{and} \quad v(I) \geq \deg(I)v(t_1) / 2 \quad \text{for } I \in S\setminus \{ H \}.
		    \]
		    Therefore,
		    \[12v(I)-\mathrm{deg}(I)v(H)\geq{0} \text{   for all   } I\in{S}.\]
		    On the other hand, computing the reduction modulo $\mfrak$ of $I_{4}/t_{1}^2$ and $I_{18}/t_{1}^{9}$, we find that 
		    \begin{align*}
    		    \overline{I_{4}/t_{1}^{2}}&=-2(\overline{\mu}_{1}^2-\overline{\mu}_{1}\overline{\mu}_{2}+\overline{\mu}_{2}^{2}),\\
    		    \overline{I_{18}/t_{1}^{9}}&=-\overline{\mu}_{1}^2\overline{\mu}_{2}^2(\overline{\mu}_{1}-2\overline{\mu}_{2})(\overline{\mu}_{1}+\overline{\mu}_{2})(2\overline{\mu}_{1}-\overline{\mu}_{2})(\overline{\mu}_{1}-\overline{\mu}_{2})^2.
		    \end{align*}
		    It is not so hard to check that, since $p\neq{3}$, the quantities $\overline{I_{4}/t_{1}^{2}}$ and $\overline{I_{18}/t_{1}^{9}}$ cannot be simultaneously zero. We thus find that 
		     \[
		    v(I_{4}) = 2 v(t_1) \quad \text{or} \quad v(I_{18}) = 9 v(t_1).
		    \]
		    Our computations give 
		    $v(\Delta / t_1^{6}) \geq 0$, so we deduce that $v(\Delta) > 4 v(t_1)$ and thus
		    \[
		        v(\Delta) - 2v(I_{4}) > 0 \quad \text{or} \quad 9v(\Delta) - 4v(I_{18}) > 0.
		    \]
		    Combining these, we obtain the statement of the theorem.

		    \subsubsection{Type III}
		    
		    Consider the universal family for Type~$\mathrm{III}.2$ in Table~\ref{tab:UniFam}. We compute the reduction of $I_{4}$ modulo $\mathfrak{m}$ to obtain
		    \[
		    \overline{I}_4 = -2.
		    \]
		    So, since $p \neq 2$, we deduce that $v(I_4) = 0$. Computing $\Delta$ and $H$, we obtain that $\Delta$ is divisible by $t_{1}^{2}t_{2}^{2}$ and $\overline{H}=0$. So we deduce that
		    \[
		         v(\Delta) - 2 v(I_4) > 0 \quad \text{and} \quad v(H) - 3 v(I_4) > 0.
		    \]

		    \subsubsection{Finishing the proof}\label{sec:SufficiencyTheorem1}
 		    Finally, we check sufficiency. Suppose that the condition in \ref{MT1:type1} is satisfied. In particular,
		    \[
		            8v(I_4)-4v(\Delta)\geq{0} \text{ and } 8v(I_{18})-18v(\Delta)\geq{0}.
		    \]
	        But these imply that
	        \[
	            v(\Delta) - 2v(I_4)\leq 0 \text{ and } 9v(\Delta) - 4v(I_{18})\leq 0,
	        \]
	        so the conditions in \ref{MT1:type2} and \ref{MT1:type3} can not be satisfied. Now suppose that the conditions in \ref{MT1:type2} are satisfied. The condition
		    \[
		        v(\Delta)- 2v(I_4) > 0 \text{ or } 9v(\Delta) - 4v(I_{18}) > 0
		    \]
		    gives
		    \[
		        8v(I_4)-\deg(I_4)v(\Delta) < 0 \text{ or } 8v(I_{18})-\deg(I_{18})v(\Delta) < 0,
		    \]
	        and therefore, the condition in \ref{MT1:type1} can not be satisfied. On the other hand, the inequality $12v(I_4) - 4v(H) \geq 0$ implies $v(H) - 3v(I_4)\leq 0$, which means that the conditions in \ref{MT1:type3} can not be satisfied either. Finally, suppose that the conditions in \ref{MT1:type3} are satisfied. Hence,
	        \[
	            8v(I_4)-\deg(I_4)v(\Delta) < 0 \text{ and } 12v(I_4)-\deg(I_4)v(H) < 0,
	        \]
	        so that the conditions in \ref{MT1:type1} and \ref{MT1:type2} can not be satisfied. This finishes the proof.

        \subsection{Proof of \texorpdfstring{\cref{thm:reductionTypes}}~}
        
        In this subsection, we prove \cref{thm:reductionTypes}. The strategy is the same as in Section~\ref{pf:Thm1}: we first calculate the invariants for each universal family and show that the given inequalities hold. This gives the necessity of the conditions. The sufficiency in this case is trivial, so this concludes the proof.   
        We start with trees of Type~$\mathrm{II}.1$ since the marking does not matter for trees of Type~$\mathrm{I}$.

        \subsubsection{Type II.1}

        We calculate $j_{2}$ and find that it is divisible by $t_{1}^2$.  Computing the reductions modulo $\mathfrak{m}$ of $j_5$ we obtain $\overline{j}_5 = 1$.
        We then deduce that $5 v(j_2) - 2 v(j_5) > 0$.
        
        \subsubsection{Type II.2}
         We calculate the reductions of $j_{2}$ and $j_{5}$ modulo $\mathfrak{m}$ and find $\overline{j}_2 = \overline{\mu}_2^2$  and $\overline{j}_5 = 1$. So we deduce that $ 5v(j_2) - 2v(j_5) = 0$.

        \subsubsection{Type III.1}
          We find that $j_{2}$ is divisible by $t_{1}^2$. Calculating the reduction of $j_{5}$ modulo $\mfrak$, we obtain $\overline{j}_5 = 1$. This implies $5 v(j_2) - 2 v(j_5) > 0$.

        \subsubsection{Type III.2}
        
        We compute the reductions of $j_{2}$ and $j_{5}$ modulo $\mfrak$ and find $\overline{j}_2 = \overline{j}_5 = 1$. This implies $5 v(j_2) - 2 v(j_5) = 0$.
        
    \subsection{Proof of \texorpdfstring{\cref{thm:edgeLengths}}~}
    
        For trees of Type~$\mathrm{I}$, there is nothing to prove.
        
        \subsubsection{Type II}
        
        We assume that the quintic $f$ has tree Type~$\mathrm{II}$ and use the  universal family for Type~$\mathrm{II}.1$ in Table~\ref{tab:UniFam}. Recall that the length $L(e_1)$ of the unique edge in the tree in this case is simply the valuation of $t_1$.

        We compute the invariants $\Delta, I_4$ and $I_{18}$, and find that $\Delta \in t_1^6{A}$, $I_4 \in t_1^2{A}$ and $I_{18} \in t_1^{9}A$, where $A$ is the universal algebra in Section~\ref{sec:UniversalFamilies}. Computing the reductions, we see
        \begin{align*}
            \overline{\Delta / t_1^6} &= \overline{\mu}_1^2  \ \overline{\mu}_2^2 \   (\overline{\mu}_1 - \overline{\mu}_2)^2, \\
            \overline{I_4/t_1^2} &= -2 \ (\overline{\mu}_{1}^{2} - \ \overline{\mu}_{1} \overline{\mu}_{2} + \ \overline{\mu}_{2}^{2}), \\
            \overline{I_{18} / t_1^9} &= - \overline{\mu}_1^2 \ \overline{\mu}_2^2 \   (\overline{\mu}_1 + \overline{\mu}_2) \ (\overline{\mu}_1 - 2 \overline{\mu}_2) \ (2 \overline{\mu}_1 - \overline{\mu}_2)  \ (\overline{\mu}_1 - \overline{\mu}_2)^2.
        \end{align*}
        Notice that, since $p \neq 3$, the two quantities $\overline{I_4/t_1^2}$ and $\overline{I_{18} / t_1^9}$ cannot vanish simultaneously. So we obtain
        \begin{align*}
            v(\Delta) = 6 v(t_1) \quad \text{and} \quad \left( \ v(I_4) = 2 v(t_1) \quad \text{or} \quad v(I_{18}) = 9 v(t_1) \ \right).
        \end{align*}
        From this, we deduce that
        \[
        v(t_1) = \frac{1}{2} \left(v(\Delta) - 2v(I_4)\right) \quad \text{or} \quad v(t_1) =  \frac{1}{3} \left(2 v(\Delta) - v(I_{18})\right),
        \]
        and the valuation of $t_1$, the length of the unique non-trivial edge, is then the maximum of these two quantities, i.e.,
        \[
            L(e_1) = \max\left(\frac{1}{2} \left(v(\Delta) - 2v(I_4)\right), \ \  \frac{1}{3} \left(2 v(\Delta) -v(I_{18})\right) \right).
        \]
        
        \subsubsection{Type III}
        
        Before we give the proof of Theorem~\ref{thm:edgeLengths} for trees of Type~$\mathrm{III}$, we shortly discuss the various formulas occurring in that theorem. For a tree of Type~$\mathrm{III}$, we can only recover the edge lengths from the quintic invariants up to a permutation of the edges, see Example~\ref{exa:TreesFiveLeaves}. A set of representatives of the edges here is given by $(e_{1},e_{2})$ with $L(e_{1})\leq{L(e_{2})}$. For trees of Type~$\mathrm{III}.2$, this symmetry continues to hold, so the formulas do not change. For trees of Type~$\mathrm{III}.1$, there is no such symmetry, meaning that we can single out the edge next to the marked point. See Section~\ref{subsubsec:TypeIII.1} for the formulas in this case.
        
       Now suppose that $f$ has tree Type~$\mathrm{III}$. Recall that the lengths $L(e_1)$ and $L(e_2)$ of the edges $e_1$ and $e_2$ are respectively $v(t_1)$ and $v(t_2)$. Computing the invariants $\Delta, I_4$ and $I_{18}$, we obtain\footnote{The bold-faced factors have positive valuation, the remaining factors have valuation $0$.}
        \begin{align*}
        \Delta = \  &{ \color{black} \bm{t_1^2 \ t_2^2}} \ \mu_1^2 \ \mu_2^2 \ (t_2 \mu_2 + 1)^2 \ (-t_1 \mu_1 + t_2 \mu_2 + 1)^2 \ (t_1 \mu_1 - 1)^2
        \end{align*}
        
        \begin{align*}
        I_{18} = \  &{ \color{black} \bm{(\mu_1 t_1 - \mu_2 t_2) (\mu_1  t_1 + \mu_2  t_2)}} (-\mu_1  t_1 + \mu_2  t_2 + 2) (\mu_2^2  t_2^2 - \mu_1  t_1 + 2  \mu_2  t_2 + 1)\\
                 \  &( \mu_2^2  t_2^2 + \mu_1  t_1 - 1)    (-2  \mu_1  \mu_2  t_1  t_2 + \mu_2^2  t_2^2 - \mu_1  t_1 + 2  \mu_2  t_2 + 1)    (\mu_1  \mu_2  t_1  t_2 - \mu_2  t_2 - 1)\\
                 \  &( \mu_1  \mu_2  t_1  t_2 + \mu_2  t_2 + 1)    (\mu_1  \mu_2  t_1  t_2 - \mu_1  t_1 + 1)    (\mu_1  \mu_2  t_1  t_2 + \mu_1  t_1 - 1) \\
                 \  &( \mu_1  \mu_2  t_1  t_2 + \mu_1  t_1 - 2  \mu_2  t_2 - 1)    (\mu_1  \mu_2  t_1  t_2 + 2  \mu_1  t_1 - \mu_2  t_2 - 1) (-\mu_1^2  t_1^2 + \mu_2  t_2 + 1)\\
                 \  &(-\mu_1^2  t_1^2 + 2  \mu_1  \mu_2  t_1  t_2 + 2  \mu_1  t_1 - \mu_2  t_2 - 1) (\mu_1^2  t_1^2 - 2  \mu_1  t_1 + \mu_2 t_2 t_2 + 1)
        \end{align*}
        and $I_4 \in A$ with $\overline{I}_4 = -2$. So we deduce that
        \[
            v(\Delta) = 2 v(t_1) + 2 v(t_2) \quad \text{and} \quad v(I_4) = 0.
        \]
        If $v(t_1) < v(t_2)$, then we find $v(I_{18}) = 2 v(t_1)$ and this gives
        \[
             L(e_1) = v(t_1) = \frac{1}{2}\left(v(I_{18}) - \frac{9}{2} v(I_{4})\right).
        \]
        If, on the other hand, $v(t_1) = v(t_2)$, then we have
        \[
            L(e_1) = v(t_1) = \frac{1}{4} \left( v(\Delta) - 2 v(I_4) \right).
        \]
        Therefore, in both cases, we obtain
        \[
            L(e_1) = v(t_1) = \min \left(\frac{1}{2}\left(v(I_{18}) - \frac{9}{2} v(I_{4})\right),  \frac{1}{4} \left( v(\Delta) - 2 v(I_4) \right) \right).
        \]
        The length of the second edge is $v(t_2)$ and can be computed using $\Delta$ as
        \[
            L(e_2) = v(t_2) = \frac{1}{2} (v(\Delta) - 2v(I_4)) -  v(t_1).
        \]
        Hence we deduce that
        \begin{align*}
            L(e_1) &= \min \left(\frac{1}{2}\left(v(I_{18}) - \frac{9}{2} v(I_{4})\right),  \frac{1}{4} \left( v(\Delta) - 2 v(I_4) \right) \right),\\
            L(e_2) &= \frac{1}{2} (v(\Delta) - 2v(I_4)) -  L(e_1).
        \end{align*}
        
        \subsubsection{Type III.1}\label{subsubsec:TypeIII.1}
        
        Now let $(q,\ell)$ be a $(4,1)$-form with tree Type~$\mathrm{III}.1$. Let $e_{1}$ be the edge adjacent to $\infty$ and $e_2$ the second edge in the tree Type~$\mathrm{III}.1$, see \cref{fig:fourOneTreeTypes}. Using the universal family for Type~$\mathrm{III}.1$ and computing the invariants we find $j_2 \in t_1^2 A^{\times}$, $j_5 = 1$, $\Delta \in t_1^6 t_2^2 A$ and $I_4 \in t_1^2 A^{\times}$ and 
        \[ 
            \overline{j_2 / t_1^2} = \overline{\mu}_1^2 , \quad \overline{\Delta / (t_1^6 t_2^2)} = \overline{\mu}_1^4 \overline{\mu}_2^2 , \quad \text{and} \quad \overline{I_4 / t_1^2} = - 2 \overline{\mu}_1^2.
        \] 
        So we deduce that 
        \[ 
            L(e_1) = v(t_1) = \frac{1}{10}(5 v(j_2) - 2 v(j_5)), 
        \] 
        and 
        \[ 
            L(e_2) = v(t_2) = \frac{1}{2} \left( v(\Delta) - 2 v(I_{4}) \right) - L(e_{1}).
        \]

\subsection{Proof of \texorpdfstring{\cref{cor:reductionTypesOfPicardCurves}}~}

    We now recall from \cite{helminckInvariantsTreesSuper} how the reduction type of a Picard curve $y^3\ell(x,z)=q(x,z)$ can be recovered from the $(4,1)$-marked tree of $(q,\ell)$. We refer the reader to \cite[Section~1.2]{helminckInvariantsTreesSuper} for the definition of the reduction type of a curve. We note here that our assumption that $K$ is algebraically closed is not restrictive. Namely, if we are interested in the reduction type of a curve over a complete discretely valued field $K$, then its reduction type is completely determined by the reduction type of the base change over $\overline{K}$, see \cite[Remark 3.6]{helminckInvariantsTreesSuper}. Finally, we note that the notion of an edge length used here is the same as the notion of \emph{thickness} of the nodal point corresponding to the edge in question, which is used in other sources.

    \medskip
    
    Let $X$ be a Picard curve over $K$. The branch locus $B$ of the covering 
    \begin{equation*}
      X\to\mathbb{P}^{1}  
    \end{equation*}
    given by $[x:y:z] \mapsto [x:z]$ is the zero locus of $q\cdot\ell$. The minimal skeleton of the marked curve $(\mathbb{P}^{1,\an},B)$ is then the $(4,1)$-marked tree of $(q,\ell)$. By applying a projective transformation, we can assume that the zero of $\ell$ is $\infty$.
    
    For tame coverings, we have that the inverse image of a skeleton is a skeleton (see \cite[Theorem~3.1]{helminckInvariantsTreesSuper} or \cite[Theorem~1.1]{HELMINCK2023skeletal}), so we obtain a map 
    \begin{equation*}
        \Sigma' \to \Sigma,
    \end{equation*}
    where $\Sigma$ is the $(4,1)$-marked tree and $\Sigma'$ is its inverse image under the morphism of Berkovich analytifications $X^{\mathrm{an}}\to\mathbb{P}^{1,\mathrm{an}}$. Consider the dehomogenized polynomial $q=q(x,1)$. The criteria in \cite[Section~3.1]{helminckInvariantsTreesSuper} allow us to reconstruct $\Sigma'$ explicitly in terms of the piecewise-linear function $ - \mathrm{log}|q| $ on $\mathbb{P}^{1, \mathrm{an}} \backslash{B} $. This function can in turn be obtained from potential theory. We then find that over an edge in $\Sigma$, there are three edges if and only if the slope of $ - \mathrm{log}|q| $ is divisible by three. If it is not divisible by three, then the length of an edge has to be divided by three. That is, the expansion factor $d_{e'/e}$ in this case is three.  This data is enough to determine the skeleton for Picard curves, as the weights of the vertices are determined by the Riemann--Hurwitz conditions. The resulting graphs can be found in \cref{fig:ReductionTypesOfPicardCurves}.

        \begin{proof}[Proof of \cref{cor:reductionTypesOfPicardCurves}]
        By \cref{thm:reductionTypes}, the $(4,1)$-marked tree type of $(q,\ell)$ is determined by the tropical invariants. The edge lengths of the $(4,1)$-marked tree type are then given by \cref{thm:edgeLengths}. To obtain the edge lengths for the curve $X$, we use the formula
        \[
            d_{e'/e}L(e')=L(e),
        \]
        where $d_{e'/e}$ is the expansion factor, see \cite[Definition 2.4, Theorem 4.23]{ABBR15}.
        \end{proof}

\end{document}